\documentclass[11pt]{article}
\usepackage{amsmath}
\usepackage{amssymb}
\usepackage{graphicx}
\usepackage{epsfig}

\begin{document}

\newtheorem{Thm}{Theorem}[section]
\newtheorem{Cor}[Thm]{Corollary}
\newtheorem{Lem}[Thm]{Lemma}
\newtheorem{Prop}[Thm]{Proposition}
\newtheorem{Def}[Thm]{Definition}
\newtheorem{rem}[Thm]{Remark}

\newtheorem{Assump}[Thm]{Assumption}
 \baselineskip 0.21in

\title{\textbf{Analysis of Optimal Thresholding Algorithms for Compressed Sensing\thanks{The work was founded by the Natural Science Foundation of China (NSFC) under the  grants 12071307, 11771003, 61571384 and 61731018. }}}


\author{Yun-Bin Zhao\thanks{Shenzhen Research Institute of Big Data, Chinese University of Hong Kong, Shenzhen,  Guangdong, China. On leave from the  School of Mathematics, University of Birmingham, Edgbaston, Birmingham B15 2TT,  United Kingdom (e-mail: {\tt
y.zhao.2@bham.ac.uk}).} ~ and  ~ Zhi-Quan Luo \thanks{
 Shenzhen Research Institute of Big Data, Chinese University of Hong Kong, Shenzhen,  Guangdong, China
  (e-mail: {\tt luozq@cuhk.edu.cn}).}}

\date{ }

\maketitle

 \textbf{ Abstract.} The optimal $k$-thresholding (OT) and optimal $k$-thresholding pursuit (OTP) are newly introduced frameworks of thresholding techniques for compressed sensing and signal approximation. Such frameworks motivate the practical and efficient algorithms called relaxed optimal $k$-thresholding ($\textrm{ROT}\omega$) and relaxed optimal $k$-thresholding pursuit ($\textrm{ROTP}\omega$) which are developed through the tightest convex relaxations of OT and OTP, where $\omega$ is a prescribed integer number.  The preliminary numerical results demonstrated in \cite{Z19} indicate that these approaches can stably reconstruct signals with a wide range of sparsity levels. However, the guaranteed performance of these algorithms with parameter $ \omega \geq 2 $ has not yet established in \cite{Z19}. The purpose of this paper is to show the guaranteed performance of OT and OTP in terms of the  restricted isometry property (RIP) of nearly optimal order for the sensing matrix governing the $k$-sparse signal recovery, and to establish the first guaranteed performance result for $\textrm{ROT}\omega$ and $\textrm{ROTP}\omega$ with $ \omega\geq 2.$  In the meantime, we provide a numerical comparison between ROTP$\omega$ and  several existing thresholding methods.  \\



\section{Introduction}
 In signal processing, one is often interested in reconstructing a signal from the measurements acquired for the signal. When the signal is sparse or can be sparsely approximated,   it is possible to reconstruct the signal from far fewer measurements than the signal length  (see, e.g., \cite{CT05, D06, E10, EK12, FR13}).  More practically, one may reconstruct the most significant information of the signal (which can be interpreted as a few  largest absolute coefficients  of the signal on its redundant bases).  This amounts to solving the following minimization problem with a sparsity constraint:
\begin{equation} \label{L0} \min_{z} \{\|Az-y\|_2^2:  ~ \|z\|_0 \leq k\},\end{equation}
where $A $ is an $m \times n$ sensing  matrix with $m<n,$   $ y: =Ax \in \mathbb{R}^m $ are the measurements of the target signal $x \in \mathbb{R}^n, $    $k$ is a prescribed integer number reflecting the interested sparsity level, and   $\|z\|_0$ is called the `$\ell_0$-norm'   counting the number of nonzero entries of  $z \in \mathbb{R}^n.$
 The model (\ref{L0})  is one of the essential models for the development of theory and algorithms for compressed sensing (see, e.g., \cite {E10, EK12, FR13, NNW17}), and it also arises in other scenarios such as the subset selection  \cite{M02, BKM16},  low-rank matrix recovery \cite{CP11, CZ13, DR16, FS19},  sparse optimization and optimal control \cite{BE13b, Z18, LB19, W19}.

Thresholding  is one of the techniques that can be used to possibly solve the problem (\ref{L0}), and it was first introduced   by Donoho and Johnstone \cite{DJ94} for signal denoising problems (see also Donoho \cite{D95}). The earlier work using this technique    can also be  found  in general areas of signal processing  \cite{FN03,KR03, SNM03} and in specific areas of compressed sensing  \cite{HGT06, BD08, BT09,BD09, BD10}.  The thresholding algorithms can  be  grouped into  soft thresholding and  hard thresholding depending on the thresholding operators.  The soft ones are usually developed from a necessary optimality condition of certain optimization problems (see \cite{DDM04, D95, E06, HGT06, FR08, VW13}).  The hard  ones can be seen as the projected Landweber iteration \cite{L51} or can be derived from the  perspective of  minimizing certain surrogate functions related to the underlying sparse optimization problems (see, e.g., \cite{BD08, DDM04, L16}). The hard thresholding methods have widely been studied in the area of  compressed sensing and  signal approximation  \cite{BD08, BD09, BD10, F11, FR13, BTW15}.  The latest development of these methods can be found in such references as \cite{B14, BFH16, KK17, SL18, ZMRC19, Z19, SBRJ19}. Although the  problem (\ref{L0})  is  usually NP-hard \cite{N95}, it does not prohibit a fast development of various computational methods for this problem. Along with thresholding,   matching pursuits (e.g., \cite{MZ93, TG07, NT09, DM09}) and convex optimization (e.g., \cite{CDS98, CT05, CWB08, ZL12, ZK15, ZL17, Z18}) are also popular methods that have been widely studied in this area.

In the family of hard thresholding methods, the iterative hard thresholding (IHT) \cite{BD08, BD09, FR13} and the hard thresholding pursuit (HTP) \cite{F11, FR13}  possess the simplest structures that are easy to implement with a low computational cost.  Compared with $\ell_1$-minimization and other state-of-art algorithms, however,  the IHT and HTP are far from being robust and stable during the course of iterations and their numerical performances are  sensitive to the choice of iterative stepsize and the sparsity level of signals.  Thus some enhancing techniques were introduced  to help stabilize  the algorithm and improve the their performances. This includes the use of certain iterative stepsizes (e.g., \cite{GK09,  BD10, C11, BTW15}) and the Nestrov's acceleration technique \cite{BT09, C11, N13, KK17, KC14}.   As pointed out in \cite{Z19}, the major drawback of existing hard-thresholding-based algorithms is the direct use of  \emph{hard thresholding operator}, denoted by $ {\cal H}_k(\cdot) ,$ which retains  the $k$ largest  magnitudes of a vector and zeroing out the remaining entries of the vector. Performing hard thresholding on a vector to generate a feasible point to the problem (\ref{L0}) is independent of its objective function. This may cause a dramatic increase instead of the decrease of the objective value in the course of iterations.  The existing enhancing techniques might help in some situations, but none of them actually serves the purpose of overcoming the intrinsic drawback of the  operator ${\cal H}_k.$

To alleviate the inherent weakness of the hard thresholding operator,  Zhao \cite{Z19} introduced a new technique called the optimal $k$-thresholding (OT), based on which a new class of  thresholding algorithms was developed, and the preliminary empirical results in \cite{Z19} indicate that the optimal $k$-thresholding method is more stable and robust for signal recovery compared with traditional IHT and HTP.    The OT technique promotes the following idea: The thresholding should be made to reduce the value of the objective function instead of being independent of the objective; when  ${\cal H}_k$ is used, it  should be applied to a $k$-compressible vector (which is nearly $k$-sparse or can be approximated by a $k$-sparse vector). Such an idea is also considered in \cite{HLS20}.  The optimal $k$-thresholding selects the best $k$ components of a vector that best  fits the measurements, and it is connected directly to the reduction of the objective value.  The initial analysis   of the basic OT algorithm and the optimal $k$-thresholding pursuit (OTP) as well as  their tightest convex relaxation counterparts called ROT and ROTP have been carried out in  \cite{Z19}. It was shown in  \cite{Z19} that the restricted isometry condition $\delta_{2k}< 0.5349 $ is  sufficient to guarantee  the convergence of  OT and OTP, and that $\delta_{3k} \leq 1/5 $ is  sufficient  for both  ROT and ROTP.
 However, the convergence  of the enhanced counterparts of ROTP, referred to as  RORT2 and ROTP3 in \cite{Z19}, has not yet established. The simulations indicate  that the ROTP2 and 3 are more robust and efficient  than ROT and ROTP for signal recovery. Thus it is important to investigate the theoretical efficiency of such an advanced development of thresholding methods.
This paper is devoted exactly to such an investigation.

The first theoretical contribution of this paper is to prove some improved guaranteed performance results for OT and OTP in terms of   restricted isometry property (RIP) of nearly optimal order for the sensing matrix governing $k$-sparse or $k$-compressible signal recovery.  These results are summarized in Theorems \ref{Thm-aa} and \ref{Thm-bb} in this paper.   The main contribution  is to establish the first guaranteed performance result for the algorithms ROTP2 and ROTP3.  This is shown in the more general setting of  $\textrm{ROT}\omega$ and  $\textrm{ROTP}\omega$ algorithms which are referred to as  the optimal $k$-thresholding algorithms performing $\omega$  times of data compressions at each iteration (see Section 2 for details).  The analysis of these algorithms is carried out in noisy scenarios which include the noiseless situation  as a special case.

The paper is organized as follows.   The algorithms are described in Section 2. The analysis of OT and OTP algorithms is given in Section 3, and the analysis of   $\textrm{ROT}\omega$ and  $\textrm{ROTP}\omega$   is carried out in Section 4. The complexity of the algorithms and performance comparison with several existing ones are discussed in Section 5.  Conclusions are given in the last section.

 \emph{Notation.}  We use $\textbf{\textrm{e}} $ to denote the vector of ones and $I$   the identity matrix. All vectors are column vectors unless otherwise specified. $\mathbb{R}^n$ is the $n$-dimensional Euclidean space, and $\{0,1\}^n  $ is the set of $n$-dimensional binary vectors. $\|x\|_2,$ $\|x\|_1$ and $\|x\|_\infty$ denote the $\ell_2$-, $\ell_1$- and $\ell_\infty$-norms of the vector $x,$ respectively.    $\textrm{supp}(x) $ denotes the  support of  $x$  which is the index set $ \{i:x_ i\not =0\}. $    Given a set $S\subseteq \{1,2,\dots, n\},$
$|S|$ denotes the cardinality of $S$ and $\overline{S}=  \{1, 2,\dots, n\} \backslash S $ is the complement of $S  $ with respect to $\{1, 2,\dots, n\}.$    Given $ x\in \mathbb{R}^n$, the vector $x_S\in \mathbb{R} ^n $ is obtained by retaining the components of $x$ supported on $S$ and setting the elements outside  $S$ to be zeros.      That is, for every $ i=1, \dots, n,$  $(x_S)_i= x_i $ if $i\in S; $ otherwise, $(x_S)_i =0. $    For vectors $x$ and $ z, $    $x \otimes z$ is the Hadamard product (entry-wise product) of $x$ and $z. $
The Hadamard product of $q$ vectors $ w^{(1)} \otimes \cdots \otimes w^{(q)} $ is   written as $ \bigotimes_{j=1}^{q}  w^{(j)} . $   The vector $x$ is said to be $k$-sparse if $\|x\|_0 \leq k. $

\section{Optimal thresholding algorithms}

    Note that $A^T (y-Ax)$ is the negative gradient of the function  $\|y-Ax\|_2^2/2. $ The classic gradient method for minimizing $\|y-Ax\|_2^2/2  $  is an iterative method generating the iterates by $ x^{p+1} = x^p+ \lambda A^T (y-Ax)  $ where $ \lambda >0 $ is a stepsize. Motivated by this classic method, to generate the iterate $ x^{p+1}$ (from the current point $ x^p$) satisfying the constraint of (\ref{L0}),  the  iterative hard thresholding (IHT) method takes the following iterative scheme \cite{BD08, BD09}:
$$ x^{p+1} = {\cal H}_k \left(x^p+ \lambda   A^T(y-Ax^p) \right). $$
For simplicity, $ \lambda $ is usually set to 1. Throughout the paper, we denote by
   $u^p: = x^p+    A^T(y-Ax^p). $
The IHT provides a basis for the development of several existing algorithms including the iterative hard thresholding pursuit (HTP) in \cite{F11}, compressive sampling matching pursuit (CoSaMP) in \cite{NT09}, subspace pursuits in \cite{DM09}, and the graded hard thresholding  in \cite{B14, BFH16}.
It was pointed out in \cite{Z19} that using ${\cal H}_k$ might increase the objective value of (\ref{L0}) yielding $\|y-A{\cal H}_k(u^p)\|_2 > \| y-A x^p\|_2 ,$  unless $ u^p $ is a $k$-compressible vector. Empirical results indicate that directly applying ${\cal H}_k$ to  non-compressible  vectors  may also cause numerical oscillation or a  slow  convergence rate of the algorithms.   To overcome such a drawback,    at a given vector $u ,$ we   consider the  minimization problem
\begin{equation}   \label {compressing}   \min_{w}   \left\{ \|y- A (u  \otimes   w) \|^2_2 : ~  \sum_{i=1}^n w_i = k , ~   w \in  \{0,1\}^n \right\},
  \end{equation}
which selects the best $k$ terms of   $u,$ which admits the smallest  objective value among  all possible choices of the $k$ terms of $u. $  The total number of $k$ terms of $u$ is $ \left(\begin{array}{c} n \\ k \end{array} \right) $  which is finite, and hence the optimal solution, denoted by $ w^*,$ of the problem (\ref{compressing}) exists.    The following definition  were first introduced in \cite{Z19}.

\begin{Def}   The  $k$-sparse vector $Z^{\#}_k (u): =  u  \otimes  w^*   $ is called the optimal $k$-thresholding of $u , $  and the operator $Z^{\#}_k(\cdot)$  is called the optimal $k$-thresholding operator.
\end{Def}

A striking difference between $ Z^{\#}_k $ and $ {\cal H}_k$ lies in that performing  $Z^{\#}_k $ is directly connected to the reduction of the  objective value of  (\ref{L0}), while the operator $ {\cal H}_k $   does not involve such a mechanism to   reduce the objective.   By optimality, $Z^{\#}_k(u)$ is the best $k$ terms of $u$ in the sense that the objective at $ Z^{\#}_k (u) $ is smaller than or equal to the objective value at any other $k$  terms of $u.$ In particular,   $ \left\| y- A  Z^{\#}_k(u) \right \|_2  \leq  \| y- A  {\cal H}_k(u) \|_2 .
  $  Let $x^p$ be the current iterate.
To solve the problem (\ref{L0}),  we can use the operator $ Z^{\#}_k  $ to generate the next iterate as follows:
$$ x^{p+1}=  Z^\#_k\left( x^p+   A^T (y-Ax^p) \right), $$
 which is referred to as the \emph{optimal $k$-thresholding} (OT) algorithm in \cite{Z19}. Combining the OT with a pursuit step (projection) is called  the OTP algorithm. The pursuit step is to solve the problem (\ref{llss}) below,  which is a least-squares problem over a restricted support set.   By the definition of $Z^\#_k ,$ the two algorithms can be explicitly described as follows.

\vskip 0.08in

\textbf{OT and OTP Algorithms} \cite{Z19}. Input $(A, y, k) $ and  an initial point $x^0\in \mathbb{R}^n. $  Perform the steps below until a  stoping criterion is satisfied:
\begin{itemize}
 \item [S1]  At  $x^p$, set $u^p : = x^p +  A^T (y-A x^p) . $   Solve the problem
  \begin{equation}   \label {HT-QP}   \min_{w}   \left\{ \|y- A (u^p  \otimes   w) \|^2_2 : ~  \sum_{i=1}^n w_i= k , ~   w \in  \{0,1\}^n \right\}.
  \end{equation}
 Let $w^*$ be the solution to  this problem.

 \item[S2]   Generate the next point $x^{p+1}  $ as follows:
 \begin{itemize} \item[] For OT algorithm, let  $ x^{p+1}=  u^p  \otimes w^* .  $
 \item[] For OTP algorithm,  set $ S^{p+1} := \textrm{supp} ( u^p  \otimes w^* ) ,$ and
let $x^{p+1}$ be the solution to
 \begin{equation} \label{llss} \min_{x} \{ \|y-A x\|_2^2:  ~ \textrm{supp} (x)\subseteq S^{p+1} \}. \end{equation}

\end{itemize}

\end{itemize}
The two algorithms share the same step  S1.  The only difference lies in the second step.  In OT, the optimal $k$-thresholding of $u^p $ is directly set to be  the next iterate  $x^{p+1},$     while the OTP use the pursuit step (\ref{llss}) to chase a  point that might   be better than $ u^p  \otimes w^*  . $
A simple stopping criterion can be a prescribed number of iterations.
The OT and OTP provide a basis from which a practical and efficient algorithm can be developed. Note that the binary optimization problem (\ref{HT-QP}) is, in general,  NP-hard   \cite{CAP08, BT18}. It is natural to consider the convex relaxation  of the problem (\ref{HT-QP}), leading to the following \emph{relaxed optional $k$-thresholding} (ROT) and the \emph{relaxed optimal $k$-thresholding pursuit} (ROTP) methods proposed first in \cite{Z19}.

\vskip 0.08in

 \textbf{ROT and ROTP Algorithms}.  Input $(A, y, k)  $ and  an initial point $x^0. $ Perform the   steps below until a   stoping criterion is satisfied:

 \begin{itemize}
 \item[S1] At  $x^p$, set $u^p : = x^p +   A^T (y-A x^p) .$
 Solve the convex optimization problem
\begin{equation} \label{con-rel}    \min_{w}   \{  \|y- A (u^p  \otimes   w) \|^2_2 : ~  \sum_{i=1}^n w_i= k , ~  0\leq w\leq \textrm{\textbf{e}} \}.
\end{equation}
 Let $w^p$ be the solution to this problem.

  \item[S2] Generate $x^{p+1} $ as follows:
 \begin{itemize}
 \item [] For ROT algorithm, let  $ x^{p+1}=  {\cal H}_k (u^p  \otimes  w^p).  $
 \item [] For ROTP algorithm, set  $ x^{\#} =  {\cal H}_k (u^p  \otimes  w^p), $
 and   let $x^{p+1}$ be the solution to
 $$ \min_{x} \{ \|y-A x\|_2^2:  ~ \textrm{supp} (x)\subseteq  \textrm{supp} (x^{\#}) \}. $$
 \end{itemize}
  \end{itemize}

The problem (\ref{con-rel}) is a convex quadratic optimization problem that can be solved   efficiently by an interior-point algorithm.   As pointed out in \cite{Z19},  although the solution $w^p $ of (\ref{con-rel}) may not be exactly $k$-sparse, but it is more compressible than the original data $u^p. $ Thus the problem (\ref{con-rel}) is referred to as a `\emph{data compressing problem}'.  To make the resulting vector more compressible so that the drawback of $ {\cal H}_k$ can be dramatically reduced, we propose the following algorithm which adopts  $\omega $ times of data compression at each iteration. The algorithms are termed $\textrm{ROT}\omega$ and $\textrm{ROTP}\omega$, respectively.

\vskip 0.08in

\textbf{$\textrm{ROT}\omega$ and $\textrm{ROTP}\omega$ Algorithm.} Input $(A, y, k) . $ Give an integer number  $\omega$ and an initial point $x^0. $ Repeat the following steps until a certain stoping criterion is satisfied:

 \begin{itemize}

  \item[S1. ] At $x^p$, let $u^p: = x^p + A^T (y-A x^p).$  Set $ \vartheta   \leftarrow   u^p. $  Perform the following loops to generate the vector $ w^{(j)}, j= 1, \ldots, \omega:$

  \textbf{for} $j =1: \omega$ \textbf{do}
    \begin{equation} \label{con-opt} \min_{w}  \{  \| y- A  ( \vartheta   \otimes w) \|_2^2:   ~~ \sum_{i=1}^n w_i= k, ~
   0\leq w \leq \textrm{\textbf{e} } \}
   \end{equation}
    ~~ to obtain a solution  $w^{(j)}$  and set  $\vartheta  \leftarrow   \vartheta \otimes w^{(j)}.$

 \textbf{end }

 \item[S2. ] Let $ x^{\#}  = {\cal H}_k (u^p \otimes  w^{(1)} \otimes  \cdots \otimes w^{(\omega)}) .  $
Generate $x^{p+1} $ as follows:
 \begin{itemize}
 \item [] For $\textrm{ROT}\omega$ algorithm, let  $ x^{p+1}=  x^{\#} .  $
 \item [] For $\textrm{ROTP}\omega$ algorithm, let  $x^{p+1}$ be the solution to the problem
 $$ \min_{x} \{ \|y-A x\|_2^2:  ~ \textrm{supp} (x)\subseteq  \textrm{supp}  ( x^{\#} ) \}. $$
  \end{itemize}

  \end{itemize}
In  step S1, we perform $\omega$ times of data compression  by solving the problem (\ref{con-opt}) starting from $ u^p. $ Specifically, after $j$th compression,  the $(j+1)$th compression  is to solve the the convex quadratic optimization problem
$$  \min_{w}  \left\{  \left\| y- A  \left[  \left ( u^p \otimes w^{(1)} \otimes  \cdots \otimes w^{(j)} \right)    \otimes w \right] \right \|_2^2:   ~~ \sum_{i=1}^n w_i= k, ~
   0\leq w \leq \textrm{\textbf{e} } \right\},
$$ to which the optimal solution is denoted by $w^{(j+1)}.$
When $ \omega=1$, the above algorithms reduce to ROT and ROPT respectively. The initial convergence results for ROT and ROTP were estalbished in \cite{Z19}. However, the convergence of ROT$_\omega$ and ROTP$_\omega$ with $ \omega \geq 2$ have not yet established. Among others, the main purpose of this paper is to  establish the first convergence result for these algorithms.

\section{Theoretical performance of OT and OTP}

The initial analysis of OT and OTP in \cite{Z19}   was performed in terms of RIP of order 2k.
In this section, we further prove that the  guaranteed performance of OT and OTP can be shown in terms of the $k$th  or $(k+1)$th order RIP of the sensing matrix. This is a nearly optimal order of RIP governing the recovery of  $k$-sparse or $k$-compressible signals.   Let us first recall the definition of RIP which has been widely used in the compressed sensing literature.

\begin{Def}  \label{Def1} \emph{\cite{CT05}} Given an $m\times n$  matrix $A  $ with $m <n,$  the $q$th order restricted isometry constant   of $A,$ denoted by $\delta_q, $ is the smallest number $\delta \geq 0$ such that
$$ (1-\delta) \|x\|^2_2 \leq \|Ax\|^2_2
\leq (1+\delta)\|x\|^2_2 $$   for any $q$-sparse vector
$x\in \mathbb{R}^n. $
\end{Def}

The following properties   will be frequently used in our later analysis.

\begin{Lem}  \label{Lem-Basic} \emph{\cite{CT05, NT09, F11}}  \emph{(i)} Let $u,v \in \mathbb{R}^n $ be $s$-sparse and $t$-sparse  vectors, respectively. If $\emph{\textrm {supp}} (u) \cap \emph{\textrm{supp}} (v) = \emptyset, $ then $$ |u^T A^T A v | \leq \delta_{s+t} \|u\|_2 \|v\|_2. $$
\emph{(ii)} Let $  v \in \mathbb{R}^n $ be a vector and   $S \subset \{1,2, \dots, n\}$ be an index set. If $|S\cup \emph{supp} (v) | \leq t,$ one has
$$ \| [(I-A^TA)v]_S\|_2  \leq \delta_t \|v\|_2 .
$$
\end{Lem}
In what follows, we show that recovering a $k$-sparse (or $k$-compressible)  signal  via OT and OTP, the RIP bound $\delta_k< \gamma^* $ or $\delta_{k+1} \leq  \gamma^*  $ is very relevant, where $ \gamma^*$ is a certain positive number smaller than 1. We   distinguish two cases: $k$ is an even number or $k$ is an  odd number.

\subsection{RIP bound for $k$ being an even number} \label{even}

Assume that $k$ is an even number and denote by $\varrho =k/2.$ The following property is of independent interest.

\begin{Lem} \label{Lem-aa} Let $z$ be a $(2k)$-sparse vector. If $k$ is an even number, then
$\|Az\|_2^2 \geq (1-3 \delta_{k}) \|z \|_2^2. $
\end{Lem}

\emph{Proof.} The $(2k)$-sparse vector $z$  can be partitioned into four $\varrho$-sparse vectors with disjoint supports:
$  z = u^{(1)} + u^{(2)}+ u^{(3)}+u^{(4)}, $ where every   $u^{(i)}$ is a $\varrho$-sparse vector and
   $ {\rm supp}(u^{(i)}) \cap {\rm supp}(u^{(j)}) =\emptyset  $ for $i\not= j  . $  Clearly,
\begin{equation} \label{EE01} \| z \|_2^2 = \sum_{i=1}^4 \| u^{(i)} \|_2^2.  \end{equation}
Since $u^{(1)}+u^{(2)}$ and $ u^{(3)}+u^{(4)}$ are $k$-sparse, by the definition of the constant $\delta_k, $ we have
 \begin{equation} \label{EE02a}  \| A(u^{(1)}+u^{(2)}) \|_2^2 \geq  (1-\delta_k) \|u^{(1)}+u^{(2)}\|_2^2 = (1-\delta_k) (\|u^{(1)}\|_2^2 +\|u^{(2)}\|_2^2), \end{equation}
  \begin{equation} \label{EE02b}  \| A(u^{(3)}+u^{(4)}) \|_2^2 \geq  (1-\delta_k)  \|u^{(3)}+u^{(4)}\|_2^2  =  (1-\delta_k) (\|u^{(3)}\|_2^2 +\|u^{(4)}\|_2^2) . \end{equation}
  Note that for every $i\in \{1,2\}$ and $ j\in \{3, 4\},$ $ {\rm supp}(u^{(i)}) \cap {\rm supp}(u^{(j)}) =\emptyset  $ for $i\not= j    $ and  $|{\rm supp} (u^{(i)}) \cup {\rm supp} (u^{(j)})| \leq 2\varrho =k . $ It follows from Lemma \ref{Lem-Basic} that
 \begin{equation}\label{EE0301}  |(u^{(i)})^T A^T A u^{(j)}| \leq \delta_ k  \|u^{(i)} \|_2  \|u^{(j)}\|_2  \leq \frac{\delta_ k}{2} (\|u^{(i)} \|_2^2+  \|u^{(j)}\|_2 ^2) . \end{equation}
Thus combining (\ref{EE01})--(\ref{EE0301}) yields
\begin{align}   \label {EEFF}  \|Az\|_2 ^2   &   =        \| A(u^{(1)} + u^{(2)}) + A( u^{(3)}+u^{(4)})\|_2^2   \nonumber \\
                       &    =       \| A(u^{(1)} + u^{(2)}) \|_2^2 + \|A( u^{(3)}+u^{(4)})\|_2^2 + 2   (u^{(1)}+u^{(2)})^T A^T A (u^{(3)}+u^{(4)})  \nonumber \\
                       &    \geq       (1-\delta_k)  \sum_{i=1}^4  \|u^{(i)}\|_2^2    + 2 [ (u^{(1)} ) ^T A^T A  u^{(3)} +  (u^{(1)}) ^T A^T A  u^{(4)} +  (u^{(2)}) ^T A^T A  u^{(3)} \nonumber \\
                       &  ~~~~+  (u^{(2)}) ^T A^T A u^{(4)} ]  \nonumber \\
                        &    \geq      (1-\delta_k) \sum_{i=1}^4  \|u^{(i)}\|_2^2  -   \delta_k [ (\|  u^{(1)}\|_2^2 +  \| u^{(3)}\|_2^2) +  (\|u^{(1)} \|_2^2+ \| u^{(4)} \|_2^2)    \nonumber\\
 &     ~~~~   +  (\|u^{(2)}\|_2^2+ \|  u^{(3)}\|_2^2)  +  (\|u^{(2)}\|_2^2+ \|u^{(4)} \|_2^2) ] \nonumber\\
                        &    =     (1-\delta_k) \sum_{i=1}^4  \|u^{(i)}\|_2^2   -   2\delta_k \sum_{i=1}^4  \|u^{(i)}\|_2^2  \nonumber\\
                       &   =    (1 - 3 \delta_k )  \|z \|_2^2,
   \end{align}
   where the first inequality follows from (\ref{EE02a}) and (\ref{EE02b}), and the second inequality follows from (\ref{EE0301}), and the final  equality follows from (\ref{EE01}). \hfill  $\Box $

\vskip 0.08in

We now prove the next technical result.

   \begin{Lem} \label{Lem-b2}  Let $h$ and $z$ be two $k$-sparse vectors, and let  $\widehat{w} \in \{0,1\}^n $ be a $k$-sparse binary vector such that $ \emph{\textrm{supp}}(h) \subseteq \emph{\textrm{supp}} (\widehat{w}). $  If $k$ is an even number, then
$$\|[(I-A^TA) (h-z)] \otimes \widehat{w} \|_2  \leq \sqrt{5}\delta_{k} \|h-z \|_2. $$
\end{Lem}

\emph{Proof.} Let $h, z, \widehat{w}$ satisfy the conditions of the Lemma.
We now partition the $k$-sparse vector $\widehat{w} $ into two binary vectors $w' $ and $ w'' $, i.e., $\widehat{w} = w' +  w'',$ where both $w'$  and $w''$ are $ \varrho$-sparse binary vectors with disjoint supports.  Note that for any vector $u\in \mathbb{R}^n$, we have
\begin{equation}\label{EEEE*} \|u \otimes \widehat{w} \|_2^2 = \|u \otimes  w' \|_2^2+  \|u \otimes  w'' \|_2^2. \end{equation}
    Let  $ v^{(1)} = (h- z) \otimes  \widehat{w} $ which is a $k$-sparse vector. Note that $ (h- z) \otimes  (\textbf{\textrm{e}} -\widehat{w})$ is also a $k$-sparse vector, and thus it can be decomposed into    $  (h- z) \otimes  (\textbf{\textrm{e}} -\widehat{w}) = v^{(2)}+v^{(3)},$   where $ v^{(2)}, v^{(3)}$ are $ \varrho$-sparse vectors with disjoint supports. Then    $ h-z= v^{(1)}+ v^{(2)}+v^{(3)}.  $ That is, $ h-z$ is decomposed into three vectors with disjoint supports.
    Since $|\textrm{supp}(v^{(1)}) \cup\textrm{supp} (\widehat{w})| \leq k,$ by Lemma  \ref{Lem-Basic}, we have
    \begin{equation} \label{ccca} \|[(I-A^TA)v^{(1)} ] \otimes  \widehat{w}  \|_2  = \|[(I-A^TA)v^{(1)} ]_{  \textrm{supp}(\widehat{w} ) }  \|_2 \leq \delta_k \|v^{(1)}\|_2. \end{equation}
    Also, we note that
       \begin{align}  \label {cccb}
       \|[(I-A^TA)(v^{(2)}+v^{(3)}) ] \otimes \widehat{w} \|_2^2  & \leq     2\left(\|[(I-A^TA)v^{(2)} ] \otimes \widehat{w} \|_2^2 + \|[(I-A^TA)v^{(3)} ] \otimes \widehat{w} \|_2^2 \right)  \nonumber \\
         & =   2 (  \|[(I-A^TA)v^{(2)} ] \otimes w' \|_2^2 +  \|[(I-A^TA)v^{(2)} ] \otimes w'' \|_2^2   \nonumber \\
          &    ~~~~ + \|[(I-A^TA)v^{(3)} ] \otimes w' \|_2^2 + \|[(I-A^TA)v^{(3)} ] \otimes w'' \|_2^2
           ) \nonumber \\
        & \leq   4 \delta_k ^2  ( \|v^{(2)} \|_2^2 +  \|v^{(3)}  \|_2^2 ),
    \end{align}
    where the first inequality follows from $\|a+b\|_2^2 \leq 2 (\|a\|_2^2+\|b\|_2^2),$ the   equality  follows from (\ref{EEEE*}), and the last inequality follows from Lemma  \ref{Lem-Basic} due to the fact $|\textrm{supp}(v^{(i)}) \cup\textrm{supp} (w ' )| \leq k $ and $|\textrm{supp}(v^{(i)}) \cup\textrm{supp} (w '')| \leq k $  for $i\in\{2,3\}. $
Then using (\ref{ccca}) and (\ref{cccb}), we have
\begin{align*}      \|[(I-A^TA)(h-z)] \otimes  \widehat{w}  \|_2
& \leq     \|[(I-A^TA)  v^{(1)}  ] \otimes  \widehat{w}  \|_2  +   \|[(I-A^TA)(v^{(2)}+v^{(3)}) ] \otimes  \widehat{w}  \|_2   \nonumber  \\
& \leq    \delta_k \|v^{(1)}\|_2 +   2 \delta_k    \sqrt{ \|v^{(2)} \|_2^2 +  \|v^{(3)}  \|_2^2 }  \nonumber  \\
&\leq    \sqrt{5} \delta_k \sqrt{\|v^{(1)}\|_2^2 + \|v^{(2)} \|_2^2 +  \|v^{(3)}  \|_2^2 }  \nonumber \\
& =   \sqrt{5} \delta_k \|z-h\|_2,
 \end{align*}
 where the third inequality follows from the fact $a + 2 \sqrt{b} \leq \sqrt{5(a^2+b)}$ for any numbers $a\geq 0$ and $ b\geq 0. $ The final  equality above follows from   $ \|z-h \|_2^2 = \|v^{(1)}\|_2^2+ \|v^{(2)}\|_2^2 + \|v^{(3)}\|_3^2. $  ~~ $\Box $

\vskip 0.08in

We now show the main result for OT and OTP under the assumption of $\delta_k.$  Throughout the remainder of the paper, we use $S$ to denote the index set of the $k$ largest absolute entries of the signal $x,$ and thus $x_S$ is the best $k$-term approximation of $ x. $  For convenience, we also define
\begin{equation} \label{WWKK}  {\cal W}^{k} : =  \{ w:~  \sum_{i=1}^n w_i =k, ~  w \in \{0, 1\}^n\}.  \end{equation}

\begin{Thm}  \label{Thm-aa} Let $y:= Ax+ \nu $ be the measurements of the signal $x$, where $\nu$ are the measurement errors.    Let $\{x^p\}  $  be the sequence  generated by OT or OTP algorithm. If $k$ is an even number and if the constant $\delta_k$  of $A$ satisfies \begin{equation} \label {con1}  \delta_k \leq 91/400 =0.2275 ,\end{equation}
(in particular, if $ \delta_k \leq 9/40 =0.225$),    then
$$ \| x^{p+1}-x_S    \|_2  \leq    \delta_{k}  \sqrt{\frac{5(1+\delta_k)}{1-3\delta_{k} }}     \| x^p-x_S \|_2 + \sqrt{\frac{1+\delta_k}{1-3 \delta_k}}\|A^T \nu'\|_2+\frac{2}{\sqrt{1-3\delta_k}}\|\nu'\|_2 ,$$ where $\nu'= A x_{\overline{S}}+\nu, $  and
$\rho :=   \delta_k \sqrt{\frac{5(1+\delta_k)}{1-3\delta_{k} }}  <1   $ is guaranteed under the condition (\ref{con1}).
\end{Thm}

\emph{Proof.}  Since $x_S $ and $x^{p+1} $ are $k$-sparse vectors, by Lemma \ref{Lem-aa}, we immediately have
\begin{equation}\label{EE04} \|A(x_S-x^{k+1})\|_2     \geq \sqrt{ 1-3 \delta_k  }  \| x_S-x^{p+1} \|_2.  \end{equation}
Let $\widehat{w} \in {\cal W}^{k}$    be a $k$-sparse binary vector such that $ \textrm{supp} (x_S) \subseteq   \textrm{supp} (\widehat{w}).  $  Then       \begin{equation} \label{eEeE}   x_S= x_S \otimes \widehat{w} . \end{equation}  Note that  $ u^p= x^p + A^T (y-Ax^p)  $ and $ y =Ax+\nu    =A x_S +\nu', $ where $ \nu'=Ax_{\overline{S}} +\nu . $   We immediately have that
\begin{equation} \label{WWZZ} x_S-u^p= x_S-x^p - A^T (y-Ax^p) = (I-A^TA) (x_S-x^p) - A^T \nu'.  \end{equation}
Since the vector  $(x_S- u^p) \otimes  \widehat{w} $ is a $k$-sparse vector,  we have
  \begin{equation} \label{EEE01} \|A [(x_S- u^p) \otimes  \widehat{w}]\|_2 \leq  \sqrt{1+\delta_k}   \|(x_S-u^p) \otimes  \widehat{w} \|_2.  \end{equation}
   Since $x^p$ is $k$-sparse, by (\ref{WWZZ}) and Lemma \ref{Lem-b2} (applying to the $k$-sparse vector $x_S$ and $x^p$), one has
\begin{align}  \label{EEE02}  \|(x_S- u^p) \otimes  \widehat{w}  \|_2     & \leq    \|[(I-A^TA)(x_S-x^p)] \otimes  \widehat{w}  \|_2  + \| (A^T \nu')\otimes \widehat{w} \|_2  \nonumber \\
 &  \leq  \sqrt{5} \delta_k \|x^p-x_S\|_2 + \| A^T \nu'\|_2.
 \end{align}
  By (\ref{eEeE}), (\ref{EEE01}) and (\ref{EEE02}) and using $ y= Ax_S +\nu',$  we have
\begin{eqnarray}  \label{EE033} \|y-A(u^p \otimes  \widehat{w}) \|_2   & = &  \|A[x_S- u^p \otimes  \widehat{w} ] + \nu' \|_2   \nonumber \\
 & = &  \|A [ (x_S-u^p)\otimes  \widehat{w}  ] + \nu' \|_2  \nonumber\\
&    \leq  &   \|A [ (x_S-u^p)\otimes  \widehat{w} ]    \|_2  + \| \nu'\|_2  \nonumber\\
&   \leq   &  \sqrt{1+\delta_k}   \|(x_S-u^p) \otimes  \widehat{w} \|_2  + \| \nu'\|_2 \nonumber\\
& \leq  &  \sqrt{1+\delta_k} \left(  \sqrt{5}\delta_{k}\right )  \|x_S-x^p\|_2 + \sqrt{1+\delta_k} \|A^T \nu'\|_2 + \| \nu'\|_2.
\end{eqnarray}
   For OT algorithm, $x^{p+1} = u^p \otimes  w^*,$ where $w^*$  is the minimizer of the problem (\ref{HT-QP}). Thus  $\|y-Ax^{p+1} \|_2 = \|y-A(u^p \otimes  w^*) \|_2.    $  For OTP algorithm,   the  iterate $x^{p+1}$ is obtained by solving the problem    $$ \min \{\|y-A z \|_2^2:   ~ \textrm{supp} (z) \subseteq  \textrm{supp} (u^p \otimes w^*)   \}, $$      which implies that   $\|y-Ax^{p+1} \|_2 \leq \|y-A(u^p \otimes  w^*) \|_2 . $
 Therefore, by optimality, the sequence $\{x^p\} $ generated by OT or OTP satisfies
\begin{equation}  \label{IINNQQ}   \|y-Ax^{p+1} \|_2 \leq \|y-A(u^p \otimes  w^*) \|_2 \leq  \|y-A(u^p \otimes  w ) \|_2  ~ \textrm{    for any } w \in {\cal W}^k. \end{equation}
In particular, since $ \widehat{w} \in {\cal W}^k, $ we have
\begin{equation} \label{equation-d}    \|y-A x^{p+1}\|_2    \leq   \|y-A(u^p \otimes  \widehat{w} ) \|_2.  \end{equation}
 By the triangle inequality, $$\|y-A x^{p+1}\|_2  = \|A(x_S-x^{p+1}) +\nu'\|_2 \geq \|A(x_S-x^{p+1})\|_2- \|\nu'\|_2. $$ which together with  (\ref{equation-d}) implies that
 \begin{equation} \label{SSDD}  \|A(x_S-x^{p+1})\|_2   \leq  \|y-A(u^p \otimes  \widehat{w} ) \|_2 + \|\nu'\|_2.  \end{equation}
  Merging (\ref{EE04}), (\ref{EE033}) and  (\ref{SSDD}) yields
 \begin{align*}  \| x_S -x^{p+1}\|_2  & \leq    \frac{1}{ \sqrt{1-3\delta_{k}  } }  (\|y-A(u^p \otimes \widehat{w} ) \|_2  + \|\nu'\|_2)  \\
  & \leq    \delta_{k}  \sqrt{\frac{5(1+\delta_k)}{1-3\delta_{k} }}     \|x_S-x^p\|_2 + \sqrt{\frac{1+\delta_k}{1-3 \delta_k}}\|A^T \nu'\|_2  + \frac{2}{\sqrt{1-3\delta_k}}\|\nu'\|_2 .
 \end{align*}
 Define $  \psi (\gamma) := 5 \gamma^3+ 5 \gamma^2+3 \gamma. $
 Clearly, $\rho: =  \delta_{k}  \sqrt{\frac{5(1+\delta_k)}{1-3\delta_{k} }} <1  $  is equivalent to   $  \psi ( \delta_k) = 5  \delta_{k}^3+ 5 \delta_{k}^2+3 \delta_{k} < 1. $ To ensure this inequality, it is sufficient to require that $\delta_k <\gamma^*,$ where  $\gamma^*$  is the real root of the univariate equation  $ \psi(\gamma) = 1 $ in the interval $[0,1].$   It is easy to verify that $ 1> \gamma^*>  91/400= 0.2275  $ and $ \psi (\gamma)$ is strictly increasing in $[0,\gamma^*].$     Thus $ \psi(\delta_{k}) < 1$ is guaranteed  if $\delta_k \leq 91/400 =0.2275$.   In particular, this is guaranteed if $\delta_k \leq 9/40 =0.225.$   \hfill    $\Box $

\subsection{RIP bound for $k$ being an odd number}

We now consider the case when the sparsity level  $k$ is an odd number, i.e.,  $k= 2 \varrho + 1 .  $
  The following lemma is similar to Lemma \ref{Lem-aa}.

\begin{Lem} \label{Lem-bb} Let $z$ be a  given $(2k)$-sparse vector, where $k= 2 \varrho + 1$ is an odd integer number.   Then
$\|Az\|_2^2 \geq (1- \delta_{k+1}- 2 \delta_{k}) \|z \|_2^2. $
\end{Lem}

\emph{Proof.} When $k= 2 \varrho + 1, $ the ($2k$)-sparse vector  $z $ can be partitioned into the following four sparse vectors with disjoint supports:
$   z = u^{(1)} + u^{(2)}+ u^{(3)}+u^{(4)}, $ where    $u^{(1)}$ and $ u^{(2)}$ are $\varrho$-sparse, $u^{(3)}$ and $  u^{(4)}$ are $(\varrho+1)$-sparse and
   $  \textrm{supp} (u^{(i)}) \cap   \textrm{supp} (u^{(j)}) =\emptyset  $ for $i\not= j  . $  Clearly,
$  \|z \|_2^2 = \sum_{i=1}^4 \| u^{(i)}\|_2^2.  $
The  two inequalities below follows immediately from   Definition   \ref{Def1}:
\begin{equation} \label{I01} \| A(u^{(1)}+u^{(2)}) \|_2^2 \geq  (1-\delta_{2\varrho}) \|u^{(1)}+u^{(2)}\|_2^2 = (1-\delta_{k-1}) (\|u^{(1)}\|_2^2 +\|u^{(2)}\|_2^2 ) ,  \end{equation}
\begin{equation} \label{I02}   \| A(u^{(3)}+u^{(4)}) \|_2^2 \geq  (1-\delta_{2(\varrho+1)}) \|u^{(3)}+u^{(4)}\|_2^2 = (1-\delta_{k+1}) (\|u^{(3)}\|_2^2 +\|u^{(4)}\|_2^2. \end{equation}
For every $i\in \{1,2\}$ and $ j\in \{3, 4\},$ one has  $|  \textrm{supp} (u^{(i)}) \cup \textrm{supp} (u^{(j)})| \leq \varrho + (\varrho +1)=k ,$ by Lemma \ref{Lem-Basic},  we have
\begin{equation} \label{I03}   |(u^{(i)})^T A^T A u^{(j)}| \leq \delta_ k  \|u^{(i)}\|_2  \|u^{(j)} \|_2  \leq (\delta_ k/2) (\|u^{(i)} \|_2^2+  \|u^{(j)}\|_2 ^2) . \end{equation}
By (\ref{I01})--(\ref{I03}) and a  similar proof to (\ref{EEFF}), we have
\begin{align*}    \|Az \|_2 ^2
                       &    =         \| A(u^{(1)} + u^{(2)}) \|_2^2 + \|A( u^{(3)}+u^{(4)})\|_2^2 + 2   (u^{(1)}+u^{(2)})^T A^T A (u^{(3)}+u^{(4)}) \\
                        &   \geq     (1-\delta_{k-1})  (\|u^{(1)}\|_2^2 + \|u^{(2)} \|_2^2) +  (1-\delta_{k+1}) (\|u^{(3)}\|_2^2 + \|u^{(4)} \|_2^2) -   \delta_k [(\|  u^{(1)}\|_2^2 +  \| u^{(3)}\|_2^2) \\
                       &     ~~~    +  (\|u^{(1)} \|_2^2+ \| u^{(4)} \|_2^2) +  (\|u^{(2)}\|_2^2+ \|  u^{(3)}\|_2^2) +  (\|u^{(2)}\|_2^2+ \|u^{(4)} \|^2_2) ]\\
                        &    =     (1-\delta_{k+1}) \sum_{i=1}^4  \|u^{(i)}\|_2^2 +  (\delta_{k+1}-\delta_{k-1})  (\|u^{(1)}\|_2^2 + \|u^{(2)} \|_2^2)  -   2\delta_k \sum_{i=1}^4  \|u^{(i)}\|_2^2 \\
                       &    \geq      (1 - \delta_{k+1}-2 \delta_k )  \| z \|_2^2,
   \end{align*}
where the last inequality follows from the fact $\delta_{k-1}\leq \delta_{k+1}.$   ~ $ \Box $

\vskip 0.08in

 We now prove the main result for OT and OTP algorithms when $ k$ is an odd number.

 \begin{Thm}  \label{Thm-bb} Let $y: = Ax+\nu $ be the measurements of the   signal $x, $  where  $\nu$ are the measurement errors. Let the sequence $\{x^p\}  $ be generated by the algorithm OT or OTP.   If $k$ is  an odd number, and if the   constant  $\delta_{k+1} $ of $A$ satisfies that
  \begin{equation} \label{ASSUMP}  \delta_{k+1} \leq  91/400= 0.2275   ,\end{equation} (in particular,   if $ \delta_{k+1} \leq 9/40 =0.225  $)
then
$$ \| x^{p+1}-x_S  \|_2  \leq    \rho \| x^p-x_S \|_2 +    \sqrt{ \frac{ 1+\delta_k  }  {   1-  \delta_{k+1}-2 \delta_k    }  } \|A^T\nu'\|_2 +   \frac{  2  }  { \sqrt{ 1-  \delta_{k+1}-2 \delta_k  }   }  \|\nu'\|_2  ,$$ where
$$  \rho : = \delta_{k+1}  \sqrt{  \frac{ 5(1+\delta_k) }   { 1-  \delta_{k+1}-2 \delta_k  } }      <1,  $$  which is ensured under the condition (\ref{ASSUMP}).

\end{Thm}

\emph{Proof.}  Since $x_S-x^{k+1} $ is $(2k)$-sparse, by setting $z= x_S-x^{p+1}$ in Lemma \ref{Lem-bb}, we immediately obtain the following relation:
\begin{equation}\label{EE044} \|A(x_S-x^{k+1})\|_2   \geq \sqrt{ 1-\delta_{k+1}-2 \delta_k  }  \| x_S-x^{p+1} \|_2.  \end{equation}
Similar to the proof of Theorem \ref{Thm-aa}, we still let
  $\widehat{w} \in {\cal W}^k   $ (which is defined by (\ref{WWKK})) be a  $k$-sparse binary vector such that $ \textrm{supp} (x_S) \subseteq   \textrm{supp} (\widehat{w})   $  and thus  $x_S =  x_S \otimes   \widehat{w}  .$   Since   $(x_S- u^p) \otimes  \widehat{w} $ is  $k$-sparse,  one has
  \begin{equation} \label{EEE1} \|A [(x_S- u^p) \otimes  \widehat{w}]\|_2 \leq  \sqrt{1+\delta_k}   \|(x_S-u^p) \otimes  \widehat{w} \|_2 .  \end{equation}
The vector $\widehat{w} $ can be partitioned as $\widehat{w} = w' +  w'', $ where $w'$ is a $ \varrho$-sparse binary vector and  $w''$ is a $ (\varrho+1)$-sparse binary vector and the supports of  $w' $ and $ w''$ are disjoint.
Partition the $(2k)$-sparse vector $x_S-x^p$ into three vectors $\eta^{(1)}, \eta^{(2)}$ and $\eta^{(3)}$  with disjoint supports    such that  $$ \eta^{(1)}= (x_S- x^p) \otimes  \widehat{w}, ~~ \eta^{(2)}+\eta^{(3)}=  (x_S- x^p) \otimes  (\textbf{\textrm{e}}-\widehat{w}) ,  $$ where $ \eta^{(2)} $ and $  \eta^{(3)}$ are $ \varrho$-sparse and $ (\varrho+1)$-sparse vectors, respectively.
Note that  $$ \begin{array}{ll}     |\textrm{supp} (\eta^{(1)})  \cup \textrm{supp}( w')|  \leq   k,  &   |\textrm{supp} (\eta^{(2)})  \cup \textrm{supp}( w')| \leq  2\varrho   = k-1 ,  \\[6pt]
  |\textrm{supp} (\eta^{(3)})  \cup \textrm{supp}( w')| \leq  2\varrho +1  = k,   &  |\textrm{supp} (\eta^{(1)})  \cup \textrm{supp}( w'')| \leq   k,  \\ [6pt]
  |\textrm{supp} (\eta^{(2)})  \cup \textrm{supp}( w'')| \leq  2\varrho +1   = k ,   &   |\textrm{supp} (\eta^{(3)})  \cup \textrm{supp}( w'')| \leq  2\varrho +2  = k+1.
  \end{array}
  $$  It  follows from Lemma \ref{Lem-Basic} that
\begin{align}   \label{aaxx}           \| [(I-     A^TA)    ( \eta^{(2)}+\eta^{(3)}) ]\otimes w'  \|_2 ^2
  & \leq    2  (\| [(I-A^TA)  \eta^{(2)}  ] \otimes w'  \|_2 ^2 + \| [(I-A^TA)  \eta^{(3)}  ] \otimes w' \|_2 ^2 ) \nonumber \\[5pt]
 &  \leq   2(\delta_{k-1}^2 \|\eta^{(2)}\|_2^2 + \delta_{k}^2 \|\eta^{(3)}\|_2^2).
 \end{align}
 Similarly,
 \begin{equation} \label{aacc}  \|  [(I-A^TA)  ( \eta^{(2)}+\eta^{(3)})  ]\otimes  w'' \|_2 ^2  \leq 2 (\delta_{k}^2 \|\eta^{(2)}\|_2^2 + \delta_{k+1}^2 \|\eta^{(3)}\|_2^2). \end{equation}
Since  $ x_S-x^p= \eta^{(1)}+ \eta^{(2)}+\eta^{(3)} ,  $  by (\ref{WWZZ}), (\ref{aaxx}) and (\ref{aacc}), we have
\begin{eqnarray}  \label{IIaa}  & & \|(x_S - u^p) \otimes  \widehat{w}  \|_2   \nonumber\\
 &  &=    \|[(I-A^TA)(x_S-x^p)] \otimes  \widehat{w}   - (A^T \nu') \otimes  \widehat{w}   \|_2    \nonumber \\
 & &=     \|  [(I-A^TA) ( \eta^{(1)}+ \eta^{(2)}+\eta^{(3)}) ]\otimes  \widehat{w}  - (A^T \nu') \otimes  \widehat{w}   \|_2   \nonumber  \\
 & &\leq     \|  [(I-A^TA)  \eta^{(1)}  ]\otimes  \widehat{w}  \|_2   +    \|  [(I-A^TA) ( \eta^{(2)}+\eta^{(3)}) ]\otimes  \widehat{w}  \|_2  +\| A^T \nu'\|_2   \nonumber  \\
 & & \leq    \delta_{k}    \|  \eta^{(1)}  \|_2  +   \left( \|  [(I-A^TA)  ( \eta^{(2)}+\eta^{(3)}) ] \otimes    w'   \|_2 ^2 +  \|  [(I-A^TA)  ( \eta^{(2)}+\eta^{(3)})  ] \otimes  w''   \|_2 ^2  \right)^{1/2}   \nonumber  \\
 & & ~~~ +  \| A^T \nu'\|_2   \nonumber  \\
 &  & \leq     \delta_{k}    \|  \eta^{(1)}  \|_2  +  \left( 2(\delta_{k-1}^2 \|\eta^{(2)}\|_2^2 + \delta_{k}^2 \|\eta^{(3)}\|_2^2) + 2( \delta_{k}^2 \|\eta^{(2)}\|_2^2 + \delta_{k+1}^2 \|\eta^{(3)}\|_2^2) \right)^{1/2} +  \| A^T \nu'\|_2    \nonumber\\
 & & =    \delta_{k}    \|  \eta^{(1)}  \|_2  +  \left(2(\delta_{k-1}^2  +  \delta_{k}^2) \|\eta^{(2)}\|_2^2 + 2(\delta_{k}^2 + \delta_{k+1}^2 ) \|\eta^{(3)}\|_2^2) \right)^{1/2} +\| A^T \nu'\|_2   \nonumber \\
 & & \leq   \delta_{k+1} \left(\|  \eta^{(1)}  \|_2 + 2\sqrt{\|\eta^{(2)}\|_2^2+ \|\eta^{(3)}\|_2^2} \right) +\| A^T \nu'\|_2  ~~~~~  (\textrm{since }  \delta_{k-1} \leq \delta_k \leq \delta_{k+1}) \nonumber\\
  &  & \leq    \delta_{k+1} \sqrt{5 ( \|  \eta^{(1)}  \|_2^2 + \|\eta^{(2)}\|_2^2+ \|\eta^{(3)}\|_2^2}  +\| A^T \nu'\|_2   \nonumber\\
   &  &=   \sqrt{5}  \delta_{k+1} \|x_S- x^p\|_2  +\| A^T \nu'\|_2  ,
 \end{eqnarray}
 where the last  inequality follows from  $a+ 2 \sqrt{b} \leq \sqrt{5(a^2+b)}$ for any $ a, b\geq 0. $
 Combining (\ref{EEE1}) and (\ref{IIaa}) and noting that  $ x_S = x_S \otimes  \widehat{w}$ yields
\begin{align} \label{IIff}  \|y-A(u^p \otimes  \widehat{w}) \|_2  &  = \|A(x_S- u^p \otimes \widehat{w}) + \nu'\|_2 = \|A[(x_S- u^p) \otimes \widehat{w}] + \nu'\|_2 \nonumber  \\
    & \leq     \sqrt{1+\delta_k}   \|(x_S-u^p) \otimes  \widehat{w} \|_2  +\| \nu'\|_2   \nonumber \\
  & \leq     \delta_{k+1}   \sqrt{5(1+\delta_k)} \|x_S- x^p\|_2  + \sqrt{ 1+\delta_k} \|A^T \nu'\|_2 + \|\nu'\|_2.
\end{align}
Let $w^*$  be a minimizer of the problem (\ref{HT-QP}).    As we have shown in the proof of Theorem \ref{Thm-aa}, the sequences $\{x^p\}$ and $ \{u^p\}$ generated by   OT and OTP algorithms satisfy the inequality (\ref{IINNQQ}), which implies (\ref{equation-d}) and  (\ref{SSDD}).
Combining  (\ref{SSDD}), (\ref{EE044}) and (\ref{IIff}) yields
\begin{align*}  \|x_S- x^{p+1}  \|_2  &  \leq    \frac{1}  { \sqrt{1-  \delta_{k+1}-2 \delta_k  } }     \|A(x_S-x^{p+1})\|_2  \\
& \leq    \frac{1}  { \sqrt{1-  \delta_{k+1}-2 \delta_k  } }   (\|y-A ( u^p \otimes \widehat{w})+  \|\nu'\|_2) \\
& \leq   \rho  \|x_S -x^p\|_2 +  \sqrt{ \frac{ 1+\delta_k}   {  1-  \delta_{k+1}-2 \delta_k  }   }  \|A^T\nu'\|_2 +   \frac{ 2   }  { \sqrt{ 1-  \delta_{k+1}-2 \delta_k  }   }  \|\nu'\|_2 ,
 \end{align*}  where   $$ \rho:=  \delta_{k+1} \sqrt{  \frac{ 5 (1+\delta_k) }   { 1-  \delta_{k+1}-2 \delta_k  } }.   $$
 Clearly,  $\rho <1 $  is equivalent to the condition $ 5\delta_k \delta_{k+1}^2+5\delta_{k+1}^2+ 2\delta_k+\delta_{k+1} <1 .$  It follows from the fact $\delta_k\leq \delta_{k+1}  $   that  $  5\delta_k \delta_{k+1}^2+5\delta_{k+1}^2+ 2\delta_k+\delta_{k+1} \leq  5  \delta_{k+1}^3+5\delta_{k+1}^2+ 3\delta_{k+1} . $ Thus the condition $\rho <1 $ is  guaranteed  if  \begin{equation} \label{IIQQ}  5  \delta_{k+1}^3+5\delta_{k+1}^2+ 3\delta_{k+1}  <1.
 \end{equation}  Let $\gamma^*$  be the real root of the univariate equation $ 5\gamma^3+5 \gamma^2+ 3\gamma =1$ in the interval $[0,1].$ It is easy to check that $ \gamma^*$ is the unique real root of this polynomial inequality in $[0,1]$ and $  5\gamma^3+5 \gamma^2+ 3\gamma <1$  for any $\gamma  \leq   \gamma^* .$  It can be verified that $\gamma^*   > 91/400. $  Therefore, $\rho < 1 $ is guaranteed under the  condition $ \delta_{k+1} \leq 91/400 $ (in particular, $ \delta_{k+1} \leq  9/40$).   \hfill  $\Box$

\vskip 0.15in

To our knowledge,    the best known RIP bound for the convergence of IHT and HTP is $\delta_{3k} \leq 1/\sqrt{3} \approx 0.5773 $  (see \cite{F11, FR13}). By adopting suitable stepsizes,   the IHT and HTP may converge under the condition $\delta_{2k} <1/3 $  (see \cite{BD10, FR13}). It was  shown in \cite{Z19} that the RIP bound for the convergence of the framework of  OT and  OTP is $ \delta_{2k} < 0.5349. $   In this section,  we have shown that OT and OTP are convergent under a nearly optimal RIP bound in terms of $\delta_k $ or $\delta_{k+1} . $

\section{Guaranteed performance of $\textrm{ROT}\omega$ and $\textrm{ROTP}\omega$}

The relaxation counterparts of OT and OTP are more practical from a computational point of view.
       The purpose of this section is to establish the first convergence result for $\textrm{ROT}\omega$ and $\textrm{ROTP}\omega$ with  $\omega \geq 2$   under the RIP assumption.   To show the main result of this section, we need to show several useful technical results which are also of independent interest. Let us start with a property of the operator ${\cal H}_k. $

\begin{Lem} \label{Lem4433}  For any vector $z \in \mathbb{R}^n $ and  any $k$-sparse vector $h\in \mathbb{R}^n,  $    one has
$$\|h- {\cal H}_k(z) \|_2 \leq \| (z- h)_{  S^* \cup S} \|_2 + \|(z- h)_{S^* \backslash S}\|_2 ,  $$
where  $S= \emph{supp} (h)  $ and  $   S^* = \emph{supp} ( {\cal H}_k(z) ). $
\end{Lem}
\emph{Proof.}      For any vector $z$, we note that   $ {\cal H}_k(z) = \textrm{arg}\min_{d}  \{ \|z- d\|_2 :   ~ \| d\|_0  \leq k\},$  which implies that
$ \| z  - {\cal H}_k(z) \|_2^2 \leq \| z - d\|_2^2 $   for any    $k$-sparse vector  $d.$
 In particular, substituting the $k$-sparse vector $d =   h+ (z-h)_S ,$ where $S= \textrm{supp} (h), $ into the inequality above leads to
  $$  \| z  - {\cal H}_k(z) \|_2^2 \leq \| z - h - (z-h)_S \|_2^2 = \| (z-h)_{  \overline{S}}\|_2^2 = \| z - h\|_2^2 - \|(z-h)_S \|_2^2 .  $$
Denote by $S^*= \textrm{supp} ({\cal H}_k(z)). $ The relation above together with  $$  \| z  -  {\cal H}_k(z)   \|^2_2
 = \|z  - h   \|^2_2  +  \|  h -  {\cal H}_k(z)    \|^2_2   -2 (h- {\cal H}_k(z)   )^T ( h-z ) .
  $$
implies that
 \begin{eqnarray*}  \label{Est-01}
 \|  h -   {\cal H}_k(z) \|^2_2
  & \leq &    - \|(z -h )_S\|_2^2 +  2 (h- {\cal H}_k(z) )^T ( h-z)    \nonumber \\
                         & =  &    - \|(z -h )_S\|_2^2 + 2  [(h- {\cal H}_k(z) )_{  S^* \cup S }]^T  ( h-z ) _{  S^* \cup S } \nonumber \\
                         & \leq   &    - \|(z-h )_S\|_2^2+ 2\|h- {\cal H}_k(z) \|_2 \|( z -h )_{  S^* \cup S }\|_2.
  \end{eqnarray*}
This further implies that $ \|  h - {\cal H}_k(z)  \|_2 $ is  smaller than or equal to the largest real root of the   quadratic equation
$ \phi(\alpha) =   \alpha^2   -  2 \alpha  \|( z -h )_{ S^* \cup S}\|_2 + \|(z -h )_S\|_2^2  =0,$
to which the largest real root is given by
\begin{align*}   \alpha^*   &  =  \left(2 \| (z -h ) _{ S^* \cup S }\|_2+ \sqrt{4 \|( z -h ) _{ S^* \cup S  } \|_2^2-4\|(z-h )_S\|_2^2}  \right)/2 \\
  &  =  \| (z- h)_{ S^* \cup S }\|_2 +  \|(z-h)_{ S^*  \backslash S}\|_2.
\end{align*}
The proof is complete. \hfill $ \Box $

 \vskip 0.08in

 The next lemma describes a property of the polytope ${\cal P} = \{w: ~\sum_{i=1}^n w_i =k, 0\leq w\leq \textbf{\textrm{e}}\}.$

 \begin{Lem} \label {Lem-B2} Let $\Lambda \subseteq \{1, \dots, n\}$ be any given index set, and let $w$ be any given vector in the polytope ${\cal P}=\{w\in \mathbb{R}^n:  ~  \sum_{i=1}^n w_i   =k, ~0\leq w\leq \emph{\textbf{\textrm{e}}}\}.$  Decompose the vector $w_\Lambda$ as the sum of $\tau$-sparse vectors:
 \begin{equation}  \label{ww-decom}  w_\Lambda  = w_{\Lambda_1} +  \cdots  +  w_{\Lambda_{q-1}}+  w_{\Lambda_{q} },  \end{equation}  where    $ \Lambda_1 \cup \dots \cup  \Lambda_{q} =  \Lambda $ and
 $\Lambda_1 $  is the index set for the $\tau $ largest  elements in $ \{w_i: i \in \Lambda \},$ and $\Lambda_2$ is the index set for the second $\tau$ largest elements in $ \{w_i: i \in \Lambda \},$    and so on. $q $ is a nonnegative integer number  such that $ |\Lambda|= (q-1)  \tau + \beta $ where $0\leq \beta < \tau . $
 Then  $$ \|w_{\Lambda _1}\|_\infty +   \cdots  + \|w_{\Lambda_{q-1}} \|_\infty +  \|w_{\Lambda_{q} } \|_\infty < (\tau+ k)/\tau .  $$

 \end{Lem}

 \emph{Proof.}     Let $ \Lambda  \subseteq \{1, \dots, n\}$ and $w\in  {\cal P} $ be given. Consider the vector $w_\Lambda $ which is decomposed as (\ref{ww-decom}). For every $i=1, \dots, q-1,$  sort the components of $ w  $  supported on $ \Lambda_i,$  i.e., $\{w_j: j\in \Lambda_i\},$  into descending order, and denote such  ordered components by $\sigma^{(i)}_1 \geq \sigma^{(i)}_2 \geq \cdots \geq \sigma^{(i)}_\tau , $ and denote the ordered components of $ w   $  supported on $ \Lambda_q $ by  $\sigma^{(q)}_1 \geq \sigma^{(q)}_2 \geq \cdots \geq \sigma^{(q)}_\beta . $
   Then the components of the vector $ w  $ supported on $ \Lambda $ are sorted into descending order as follows:
       \begin{equation} \label{sort}   \underbrace{\sigma^{(1)}_1 \geq \sigma^{(1)}_2  \geq \cdots \geq \sigma^{(1)}_\tau  }  \geq   \underbrace{\sigma^{(2)}_1 \geq \sigma^{(2)}_2  \geq \cdots \geq \sigma^{(2)}_\tau  } \geq \cdots \geq  \underbrace{\sigma^{(q)}_1 \geq \sigma^{(q)}_2  \geq \cdots \geq \sigma^{(q)}_\beta  }.   \end{equation}
      Clearly, for every $i =1, \dots, q,$  $\sigma^{(i)}_1$ is the largest entries of $ w_{\Lambda _i } , $ i.e., $\sigma^{(i)}_1 =  \|w_{\Lambda _i}\|_\infty. $   For every $i =1, \dots, q-1,$     $ \sigma_{\tau} ^{(i)} $ is the smallest entry of $ w $  on the support $ \Lambda _i,$     and $ \sigma^{(q)}_{\beta}$ is the smallest component of $w $ supported on $ \Lambda_{q}. $
      Therefore, \begin{equation} \label{DDD8}  \Phi(w, \Lambda) :=  \|w_{\Lambda _1}\|_\infty +   \cdots  + \|w_{\Lambda_{q-1}} \|_\infty +  \|w_{\Lambda_q} \|_\infty   =  \sum _{ i=1}^{q}  \sigma^{(i)}_1 . \end{equation}
It is sufficient to show that $   \Phi(w, \Lambda)  <  (\tau+ k)/\tau.  $ From (\ref{sort}), for each $i, $ the largest entry of  $ w $ on the support $\Lambda_{i+1}   $  is smaller than or equal to the smallest entry of $  w$ on the support  $  \Lambda _i , $ i.e., $ \sigma^{(i)}_\tau \geq \sigma^{(i+1)}_1 $ for every $ i \in \{ 1, \dots, q-1\}. $ So we immediately see that
\begin{equation} \label{sum}  \Phi(w, \Lambda)    =   \sum _{ i=1}^{q}  \sigma^{(i)}_1  \leq  \sigma^{(1)}_1 + \sigma^{(1)}_\tau  + \sigma^{(2)}_\tau + \cdots + \sigma^{(q-1 )}_\tau  \leq 1+  \sum_{i=1}^ {q-1} \sigma^{(i)}_\tau,
\end{equation}     where the last inequality follows from $\sigma^{(1)}_1 \leq 1$ (since $ 0\leq  w  \leq \textrm{\textbf{e}}  $).   Note that for every $i=1 , \dots, q-1, $ $ \sum_{i=1}^ {q-1} \sigma^{(i)}_\tau $ is the sum of the smallest  entries of the vector $  w $ supported on $ \Lambda_i .  $  We see from (\ref{sort}) that $$ \sum_{i=1}^ {q-1} \sigma^{(i)}_\tau \leq \sum_{i=1}^ {q -1} \sigma^{(i)}_{\tau-1} \leq \cdots \leq \sum_{i=1}^ {q -1} \sigma^{(i)}_2,   $$ which together with (\ref{sum}) implies  that
 $    \Phi(w, \Lambda)   \leq  1+  \sum_{i=1}^ {q-1} \sigma^{(i)}_j $ for $ j= 2, \dots, \tau.  $
Adding up these $\tau-1 $ inequalities and  equality (\ref{DDD8}) altogether yields
\begin{eqnarray*}  \tau  \Phi(w, \Lambda)   & \leq  &  \tau-1 + \sum_{i=1}^ {q } \sigma^{(i)}_1  +  \sum_{i=1}^ {q-1 } \sigma^{(i)}_2 + \cdots + \sum_{i=1}^ {q-1 } \sigma^{(i)}_\tau   \leq     \tau-1 +  \sum_{j\in  \Lambda }   w  _j  \leq    \tau -1+ k.  \end{eqnarray*}
  Therefore $ \Phi(w, \Lambda) \leq   ( \tau + k-1)/ \tau  < (\tau+k)/\tau, $  as desired.  ~~ $\Box $

\vskip 0.07in

We now show a property of the vectors $w^{(j)}$ generated at S1 of $ \textrm{ROT}\omega$ and $ \textrm{ROTP}\omega.$

\begin{Lem} \label{Lem4444} Let $y := Ax +\nu  $ be the measurements of $x\in \mathbb{R}^n$, where $ \nu$ are the measurement errors. Let   $ \widehat{w} \in {\cal W}^{(k)}    $ (which is defined by (\ref{WWKK})) be a binary vector such   that $  \textrm{supp}  (x_S) \subseteq   \textrm{supp}  (\widehat{w}) . $ At the iterate $x^p,$     the vectors $w^{(1)},  \cdots, w^{(\omega)}$ are generated by $ \textrm{ROT}\omega$ or $ \textrm{ROTP}\omega. $ Then
\begin{equation} \label{LL4444} \|y- A[u^p\otimes (\bigotimes_{j=1}^\omega w^{(j)} ) ] \|_2 \leq \|y-  A(u^p\otimes \widehat{w} ) \|_2 + \sum_{i=1}^ {\omega-1} \|A [(u^p-x_S) \otimes  ( \bigotimes_{j=1}^{i}  w^{(j)}  ) \otimes (\emph{\textbf{\textrm{e}}}-\widehat{w}) ] \|_2 , \end{equation}
where $u^p = x^p + A^T (y-Ax^p).  $
\end{Lem}

\emph{Proof.}  Note that $ \widehat{w} \in {\cal W}^{(k)} $  satisfies   $ \textrm{supp} (x_S) \subseteq  \textrm{supp} (\widehat{w}) . $ The first inequality below follows from the optimality of $w^{(\omega)}:$
\begin{align*}       \|y- A[u^p\otimes (\bigotimes_{j=1}^\omega w^{(j)} )  ]  \|_2
 &   \leq     \|y- A[u^p\otimes (\bigotimes_{j=1}^{\omega-1}  w^{(j)} )   \otimes  \widehat{w} ] \|_2 \\
                             &      =  \|y- A[u^p\otimes (\bigotimes_{j=1}^{\omega-1}  w^{(j)} )]    -  A [ u^p\otimes (\bigotimes_{j=1}^{\omega-1}  w^{(j)} ) \otimes (\textbf{\textrm{e}} - \widehat{w})  ] \|_2 \\
                             &      \leq   \|y- A[u^p\otimes (\bigotimes_{j=1}^{\omega-1}  w^{(j)} )]\|_2+ \|A [ u^p \otimes (\bigotimes_{j=1}^{\omega-1}  w^{(j)} ) \otimes (\textbf{\textrm{e}} - \widehat{w})  ] \|_2\\
                             &      = \|y- A[u^p\otimes (\bigotimes_{j=1}^{\omega-1}  w^{(j)} )]\|_2+ \|A [ (u^p-x_S) \otimes(\bigotimes_{j=1}^{\omega-1}  w^{(j)} ) \otimes (\textbf{\textrm{e}} - \widehat{w})  ] \|_2,
\end{align*}
where the final equality follows from $x_S \otimes (\textbf{\textrm{e}} - \widehat{w}) =0$ due to  $\textrm{supp} (x_S) \subseteq  \textrm{supp} (\widehat{w}).$
Similarly, by the optimality of $w^{(\omega-1)}, \ldots, w^{(2)},$  we obtain the following inequalities for every $\ell= \omega-1, \omega-2, \ldots, 2 : $
 $$ \|y- A[u^p\otimes (\bigotimes_{j=1}^\ell w^{(j)} )  ]  \|_2 \leq    \|y- A[u^p\otimes (\bigotimes_{j=1}^{\ell-1}  w^{(j)} )]\|_2+ \|A [ (u^p-x_S) \otimes(\bigotimes_{j=1}^{\ell-1}  w^{(j)} ) \otimes (\textbf{\textrm{e}} - \widehat{w})  ] \|_2.  $$
Merging the above inequalities altogether leads to the following relation:
$$   \|y- A[u^p\otimes (\bigotimes_{j=1}^\omega w^{(j)} )  ]  \|_2 \leq  \|y-  A(u^p\otimes w^{(1)}  ) \|_2 + \sum_{i=1}^ {\omega-1} \|A [(u^p-x_S) \otimes  ( \bigotimes_{j=1}^{i}  w^{(j)}  ) \otimes (\textbf{\textrm{e}}-\widehat{w}) ] \|_2 . $$
By the optimality of $w^{(1)}$, we have $$  \|y-  A(u^p\otimes w^{(1)}  ) \|_2 \leq \|y-  A(u^p\otimes  \widehat{w}   ) \|_2. $$
Combining the last two inequalities  above yields the desired relation (\ref{LL4444}).  \hfill    $ \Box $

\vskip 0.1in

We now bound the right-hand side of (\ref{LL4444}). The idea used to show the following lemma is based on the technique of splitting and estimation of the tail which has widely been used in the compressed sensing literature (see, e.g., Cand\`es and Tao \cite{CT05}, Foucart and Rauhut \cite{FR13}, and  Rauhut and Ward \cite{RW16}).

\begin{Lem} \label{Lem4455} Under the   conditions of Lemma \ref{Lem4444}, that is, the vectors $ x_S, \widehat{w}, x^p, u^p $ and $ w^{(j)} $ ($j=1, \cdots,  \omega) $ are the same as in Lemma \ref{Lem4444}.  Then for every $ i=1, \ldots, \omega -1, $
\begin{equation} \label{theta-bound} \Theta^{(i)}:  = \|A [(u^p-x_S) \otimes  ( \bigotimes_{j=1}^{i}  w^{(j)}  ) \otimes (\emph{\textbf{\textrm{e}}}-\widehat{w}) ] \|_2 \leq 2 \delta_{3k} \sqrt{1+ \delta_k} \|x_S-x^p\|_2 + 2\sqrt{1+\delta_k} \|A^T\nu' \|_2  ,  \end{equation}
 and thus
\begin{align} \label{Bound11} \|y- A[u^p\otimes (\bigotimes_{j=1}^\omega w^{(j)} ) ] \|_2    \leq       [\delta_{2k}    &    + 2(\omega-1) \delta_{3k}   ] \sqrt{1+ \delta_k} \|x_S-x^p\|_2  + \|\nu' \|_2 \nonumber  \\
 &   +  (2\omega-1) \sqrt{1+\delta_k} \|A^T\nu' \|_2  .
 \end{align}

\end{Lem}

\emph{Proof.} The first term in (\ref{LL4444}) is easy to bound.  Note that $y=  A x_S +\nu' $ with  $ \nu'= Ax_{\overline{S}}+\nu,$ and  $ x_S  = x_S\otimes \widehat{w} .  $ By (\ref{WWZZ}) we have
\begin{align} \label{EE03}  \|y-A(u^p \otimes  \widehat{w}) \|_2   & =     \|A [ (x_S-u^p)\otimes  \widehat{w} ]   +\nu' \|_2 \nonumber  \\
&  \leq  \|A  [(x_S-x^p) \otimes  \widehat{w} ] \|_2 +\|\nu' \|_2   \nonumber\\
& \leq \sqrt{1+\delta_k} \| (x_S-x^p) \otimes  \widehat{w} ] \|_2 +\|\nu' \|_2   \nonumber\\
&   \leq      \sqrt{1+\delta_k}  \left( \|  [(I-A^TA) (x_S-x^p)] \otimes  \widehat{w}  \|_2  +  \|(A^T\nu')\otimes  \widehat{w}  \|_2\right) + \|\nu' \|_2  \nonumber\\
& \leq     \delta_{2k} \sqrt{1+\delta_k}   \|x_S-x^p\|_2 + \sqrt{1+\delta_k} \|A^T\nu' \|_2  + \|\nu' \|_2.
\end{align}
The second inequality  above follows from  the  $k$-sparsity of the vector $(x_S-x^p) \otimes  \widehat{w} ,$  and the last inequality follows from Lemma \ref{Lem-Basic} since $
|\textrm{supp} ( x_S-x^p) \cup \textrm{supp} (\widehat{w}) |\leq 2k. $ From (\ref{LL4444}), in order to show (\ref{Bound11}),
it is sufficient to show the bound  (\ref{theta-bound})  for  $ \Theta^{(i)}, $   $i= 1, \ldots, \omega -1. $
Note that $ \textrm{supp} (\textbf{\textrm{e}}-\widehat{w}) = \overline{\textrm{supp} ( \widehat{w})}.$   Then $ \Theta^{(i)} $ can be written as
$$ \Theta^{(i)}  =  \|A [(u^p-x_S)\otimes (\bigotimes_{j=1} ^{i} ) ]_{\textrm{supp} (\textbf{\textrm{e}} -\widehat{w})} \|_2  = \| A [(u^p-x_S) \otimes  (\bigotimes_{j=1} ^{i} ) ]_{\overline{\textrm{supp} (\widehat{w}}) }  \|_2 . $$
  Let  $(w^{(1)}  )_{ \overline{\textrm{supp} ( \widehat{w} ) } } $ be decomposed into $k$-sparse vectors as follows:
   $$ (w^{(1)}  )_{ \overline{\textrm{supp} ( \widehat{w} ) } }  = (w^{(1)}  )_{T_1} +  \cdots  + (w^{(1)} )_{T_{q-1}}+  (w^{(1)}  )_{T_{q} },  $$
  where  $ {T_1}$ is the index set for  the $k$ largest elements  in  the set $ \left\{(w^{(1)}  )_i: ~ i \in   \overline{\textrm{supp}  \widehat{w} ) } \right\},$    and $T_2$ is the index set for the second $k$ largest elements in  this set, and so on.  These index sets are mutually disjoint    and the cardinality $  |T_\ell|=k \textrm{ for all } \ell =1, \dots, q-1  $ and $ |T_{q}|= \kappa' <  k, $ where  $q $ and $  \kappa' $ are integer numbers.  As a result, we have   $  \overline{\textrm{supp} ( \widehat{w} ) }     = T_1 \cup T_2 \cup  \cdots  \cup T_{q }$   with cardinality $ \left |\overline{ \textrm{supp} ( \widehat{w} ) }  \right |  =  (q-1)  k + \kappa'. $
  Note that $ w^{(1)} \in  \{w: ~ \sum_{i=1}^n w_i =k, ~ 0 \leq w \leq \textrm{\textbf{e}} \}. $     Applying   Lemma   \ref{Lem-B2} to the vectors $w= w^{(1)}, \tau = k $ and $\Lambda =\overline{\textrm{supp} ( \widehat{w} ) } ,$   we immediately have that
 \begin{equation}  \label{TTT01}  \sum_{i=1}^{q} \| (w^{(1)})_{T_i}\|_\infty < 2 . \end{equation}
Define the vector $ v^{(\ell)}: =    [(u^p-x_S) \otimes  (\bigotimes_{j=1} ^{i} w^{(j)})  ]_{T_\ell},  \ell = 1, \dots, q, $ which are $ k$-sparse vectors. Then $$  [(u^p-x_S) \otimes  (\bigotimes_{j=1} ^{i} w^{(j)} )   ]_{  \overline{\textrm{supp} ( \widehat{w} ) }  } = v^{(1)}+ v^{(2)} + \cdots +  v^ {(q)}  . $$
Therefore,
\begin{equation} \label{TTT02}  \Theta^{(i)}   =   \left\| A \sum _{ \ell=1}^ {q} v^{(\ell)} \right\|_2 \leq \sum _{ \ell =1}^ {q} \|A v^{(\ell)} \|_2 \leq \sqrt{1+\delta_{k}} \sum _{ \ell=1}^{q} \|v^{(\ell)} \|_2.   \end{equation}
 We now estimate the term  $\sum _{ \ell=1}  ^{q} \|v^{(\ell)} \|_2  .$  By the structure of the algorithm,   $w^{(j)} \in  \{w: ~ \sum_{i=1}^n w_i =k, ~ 0 \leq w \leq \textrm{\textbf{e}} \}  $ for $ j=1, \ldots, \omega.$ Thus, by (\ref{WWZZ}), we have
\begin{align}  \label{4242} \|v^{(\ell)} \|_2   & =     \|[ (u^p-x_S) \otimes (\bigotimes_{j=1} ^{i}  w^{(j)})  ]_{T_\ell }\|_2  =  \|[((I-A^TA)(x_S-x^p)- A^T \nu') \otimes (\bigotimes_{j=1} ^{i} w^{(j)})  ]_{T_\ell} \|_2  \nonumber \\
& \leq  \|[((I-A^TA)(x_S-x^p)) \otimes (\bigotimes_{j=1} ^{i} w^{(j)})  ]_{T_\ell} \|_2  + \|[(A^T \nu')\otimes (\bigotimes_{j=1} ^{i} w^{(j)})  ]_{T_\ell} \|_2 \nonumber \\
& \leq   \|(\bigotimes_{j=1} ^{i} w^{(j)})   _{T_\ell}\|_\infty \left(\|[(I-A^TA)(x_S-x^p)]_{T_\ell} \|_2  +  \| (A^T \nu')_{T_\ell}\|_2 \right ) \nonumber \\
  & \leq   \|(w^{(1)}  )_{T_\ell}\|_\infty   \left (\delta_{3k} \|x_S-x^p \|_2 + \|A^T \nu'\|_2 \right),
  \end{align}  where the   inequalities above follows from  Lemma \ref{Lem-Basic} with  $ |T_\ell \cup \textrm{supp} (x_S-  x^p)| \leq 3k , $ and from the fact $ 0 \leq w^{(j)}  \leq  \textbf{\textrm{e}}   $ for $j=1, \ldots i, $ which implies that $ (\bigotimes_{j=1} ^{i} w^{(j)})   _{T_\ell} \leq (w^{(1)}  )_{T_\ell} . $
 Thus merging (\ref{TTT01}), (\ref{TTT02}) and (\ref{4242}) yields
\begin{align*}   \Theta^{(i)}    &   \leq     \sqrt{1+\delta_k}  \left[ \sum _{ \ell =1}  ^{q}  \|(w^{(1)}  )_{T_\ell }\|_\infty   \left(\delta_{3k} \|x_S-x^p \|_2+ \|A^T \nu'\|_2\right ) \right]  \\
 & \leq 2 \delta_{3k} \sqrt{1+\delta_k}  \|x_S-x^p\|_2 + 2\sqrt{1+\delta_k} \|A^T \nu'\|_2,
\end{align*}
which is exactly the estimate in (\ref{theta-bound}).
Substituting (\ref{EE03}) and the bound above   into (\ref{LL4444}) leads to the desired bound (\ref{Bound11}). \hfill  $\Box $

\vskip 0.08in

We now prove the main result of this section which implies that the   $\textrm{ROT}\omega$ and $\textrm{ROTP}\omega$  can  recover   $x_S$ if   $ \delta_{3k}  $ is smaller than a certain number in  $ (0,1)   $ and the measurements of the signal are accurate enough.

\begin{Thm}  \label{Thm-OIHT-02}
   Let $y := Ax +\nu  $ be the measurements of $x\in \mathbb{R}^n$, where $ \nu$ are the measurement errors. Let   $\omega \geq 1 $ be a given integer number.
   \begin{itemize}\item [ \emph{(i)} ] If the restricted isometry constant of $A$ satisfies $  \delta_{3k}  < \gamma (\omega),  $ where $ \gamma (\omega)   $ is the unique real root in the interval $(0,1) $ of the   univariate equation
 $    (2 \omega +1) \gamma  \sqrt{ \frac{ 1+ \gamma  } { 1-\gamma}  }  +   \gamma =1  ,  $ then the sequence $\{x^{p}\} , $ generated by $\textrm{ROT}\omega,$  approximates $x_S$ with error   $$ \|x^{p+1}-x_S \|_2 \leq  \widetilde{\rho}  \| x^p-x_S\|_2  + c_1 \|A^T \nu'\|_2+ c_2\|\nu'\|_2 ,  $$
where $ \nu'= A x_{\overline{S}} +\nu, $ and  $$ \widetilde{\rho } : =  ( \delta_{2k}      + 2\omega  \delta_{3k} ) \sqrt{ \frac{  1+ \delta_k}   {   1-\delta_{2k}   } } +  \delta_{3k} <1
 $$ which is guaranteed under the condition $  \delta_{3k}  < \gamma (\omega).  $ The constants $ c_1$ and $ c_2$  are given by  $$ c_1=   \frac{2\omega -1}{1-\alpha^*}  \sqrt{ \frac{  1+ \delta_k}   {   1-\delta_{2k}   } }  +1 , ~~ c_2 =  \frac{2}{(1-\alpha^*)\sqrt{1-\delta_{2k}}}, $$ where $ \alpha^*= \frac{2\delta_{3k} }{ 2\omega \delta_{3k} + \delta_{2k}}. $

 \item[\emph{(ii)}] If   $\delta_{3k} <  \gamma^*(\omega),$ where $\gamma^*(\omega)    $ is the unique real root  in the interval $(0,1)$ of the univariate  equation
   \begin{equation} \label{uni-equ} \frac{1}{\sqrt{1- \gamma ^2}} \left((2 \omega +1) \gamma   \sqrt{ \frac{ 1+ \gamma  } { 1- \gamma}  }  +  \gamma \right)   =1,   \end{equation}  then the sequence $\{x^{p}\}, $ generated by $\textrm{ROTP}\omega,$  approximates $x_S$ with error
 $$  \| x^{p+1}-x_S\|_2  \leq  \rho' \|x^p -x_S \|_2 +  \tau_1 \|A^T\nu'\|_2 + \tau_2 \|\nu'\|_2 ,  $$
   where  $ \nu'= A x_{\overline{S}} +\nu,  $ and the constant \begin{equation} \label{TTEE}  \rho' :  =\frac{1}{\sqrt{1-\delta_{2k}^2}} \left( ( \delta_{2k}      + 2\omega  \delta_{3k} ) \sqrt{ \frac{  1+ \delta_k}   {   1-\delta_{2k}   } } +  \delta_{3k}     \right)  < 1    \end{equation} which is guaranteed under the condition $  \delta_{3k}  < \gamma^* (\omega). $  The constants  $ \tau_1 $ and $  \tau_2$ are given by
   \begin{equation}  \label{TTEE02}  \left\{\begin{array} {l} \tau_1=  \frac{1}{\sqrt{1-\delta_{2k}^2}} \left[\frac{2\omega -1}{1-\alpha^*}  \sqrt{ \frac{  1+ \delta_k}   {   1-\delta_{2k}   } }  +1 \right] \\
    \tau_2= \frac{1}{\sqrt{1-\delta_{2k}^2}}\left[\frac{2}{(1-\alpha^*)\sqrt{1-\delta_{2k}}} +\frac{1}{1-\delta_{2k}}\right], \end{array}  \right.   \end{equation}
   where $ \alpha^*$ is the same constant in (i).
   \end{itemize}
\end{Thm}

  \emph{Proof.}   At the current iterate $x^p$,      both $\textrm{ROT}\omega$ and $\textrm{ROTP}\omega$    generate  the   vectors $w^{(1)}, \ldots, w^{(\omega)} \in \{w: ~ \sum_{i=1}^n w_i =k, ~ 0 \leq w \leq \textrm{\textbf{e}}\} $ by solving the  convex optimization problems  (\ref{con-opt}), where $u^p= x^p+ A^T (y-Ax^p) . $     Denote by
$$x^{\#} =  {\cal H}_k(u^p \otimes (\bigotimes_{j=1}^{\omega} w^{(j)}) ), ~    X =\textrm{supp} (x^{\#}) .$$
Since  $x_S$ is a $k$-sparse vector  and $y= Ax_S+\nu' , $ where $ \nu' = A x_{\overline{S}} +\nu.  $   By Lemma \ref{Lem4433}, we have
 \begin{equation} \label{Est-111} \|  x_S -x^{\#} \|_2 \leq   \|( u^p \otimes (\bigotimes_{j=1}^{\omega} w^{(j)}) -x_S ) _{  X \cup S} \|_2 +  \|[ u^p \otimes (\bigotimes_{j=1}^{\omega} w^{(j)}) -x _S] _{ X \backslash S} \|_2  .  \end{equation}
  The second term of the right-hand side of (\ref{Est-111}) is easy to bound.  Using the fact $ (x_S)_{ X \backslash S} =0   $  and (\ref{WWZZ}), we have
\begin{eqnarray*}      \|[ u^p \otimes (\bigotimes_{j=1}^{\omega} w^{(j)}) -x_S ] _{ X \backslash S}  \|_2 &  = &    \|[ u^p \otimes (\bigotimes_{j=1}^{\omega} w^{(j)}) ] _{ X \backslash S}   \|_2   =  \|[( u^p-x_S) \otimes (\bigotimes_{j=1}^{\omega} w^{(j)}) ] _{ X \backslash S}\|_2  \nonumber \\
                         &  = &   \| [((I-A^TA) (x^p-x_S)- A^T \nu')\otimes (\bigotimes_{j=1}^{\omega} w^{(j)})]_{ X \backslash S}  \|_2  \nonumber \\
     & \leq   &  \left \| [(I-A^TA) (x^p-x_S)]_{ X \backslash S}\right\|_2 + \|A^T \nu'\|_2 \nonumber \\
     & \leq  &  \delta_{3k} \|x^p-x_S\|_2 + \|A^T \nu'\|_2 ,
\end{eqnarray*}
where the  inequalities above follow  from $ 0\leq w^{(j)}  \leq \textrm{\textbf{e}} $ for $j=1, \ldots, \omega $  and from Lemma \ref{Lem-Basic} with the fact  $|\textrm{supp} ( x^p-x) \cup (X \backslash S) |\leq 3k.$  Substituting the bound above into (\ref{Est-111}) yields
\begin{equation} \label{Bd-01} \|  x_S -x^{\#} \|_2 \leq   \|[ u^p \otimes (\bigotimes_{j=1}^{\omega} w^{(j)}) -x _S] _{  X \cup S} \|_2 +  \delta_{3k} \|x^p-x_S\|_2 + \|A^T \nu'\|_2.  \end{equation}
 We now bound the first term of the right-hand side of (\ref{Bd-01}). Let $\alpha \in (0,1)$ be any given number. Define
 $$ \Theta^* :=  \| A [ ( u^p \otimes  (\bigotimes_{j=1}^{\omega} w^{(j)}) -x _S)_{\overline{X\cup S}}  ] \|_2 . $$
 There are only two cases.

 \textbf{\emph{Case 1.}}   $ \Theta^* \leq \alpha \|A [  (u^p \otimes  (\bigotimes_{j=1}^{\omega} w^{(j)}) -x_S )_{ X \cup  S} ] \|_2. $
 In this case, since $y=Ax_S +\nu', $ we  have
 \begin{align*}        \|y-     &   A[u^p \otimes  (\bigotimes_{j=1}^{\omega} w^{(j)})   ] \|_2  \\
       & =      \|A[ u^p \otimes  (\bigotimes_{j=1}^{\omega} w^{(j)})  -x_S] -\nu' \|_2      \\
                               &       = \|A [  (u^p \otimes  (\bigotimes_{j=1}^{\omega} w^{(j)})-x_S )_{X\cup S }] + A [ ( u^p \otimes  (\bigotimes_{j=1}^{\omega} w^{(j)}) -x_S )_{\overline{X\cup S}} ]  -\nu'\|_2
       \end{align*}
       By the triangle inequality, we have
 \begin{align}   \label{DDD2E}
           \|y-     &   A[u^p \otimes  (\bigotimes_{j=1}^{\omega} w^{(j)})   ] \|_2  \nonumber  \\
         &   \geq  \|A [  (u^p \otimes  (\bigotimes_{j=1}^{\omega} w^{(j)}) -x_S )_{ X \cup  S} ] \|_2-  \| A [  (u^p \otimes  (\bigotimes_{j=1}^{\omega} w^{(j)})-x_S )_{\overline{X\cup S}} ]  \|_2  -\|\nu'\|_2\nonumber  \\
         &   =  \|A [  (u^p \otimes  (\bigotimes_{j=1}^{\omega} w^{(j)}) -x_S )_{ X \cup  S} ] \|_2-  \Theta^*   -\|\nu'\|_2 \nonumber \\
         & \geq (1-\alpha) \|A [  (u^p \otimes  (\bigotimes_{j=1}^{\omega} w^{(j)}) -x_S )_{ X \cup  S} ] \|_2  -\|\nu'\|_2\nonumber   \\
            &      \geq (1-\alpha) \sqrt{1-\delta_{2k}}  \|  [u^p \otimes  (\bigotimes_{j=1}^{\omega} w^{(j)})-x _S]_{X\cup S  } \|_2  -\|\nu'\|_2,
\end{align} where the last inequality follows from
 \begin{equation} \label{EEAA} \|A [  (u^p \otimes  (\bigotimes_{j=1}^{\omega} w^{(j)}) -x_S )_{ X \cup  S} ] \|_2 \geq \sqrt{1-\delta_{2k}} \|[u^p \otimes  (\bigotimes_{j=1}^{\omega} w^{(j)}) -x _S]_{ X \cup  S}\|_2 ,
 \end{equation}  which is implied from
 Definition \ref{Def1} and the fact $|X \cup S| \leq 2k. $
The inequality (\ref{DDD2E}), together with Lemma \ref{Lem4455}, implies that
\begin{align} \label{Est-02}     \| [u^p \otimes  (\bigotimes_{j=1}^{\omega} w^{(j)}) -x_S ]_{ X \cup  S }   \|_2
    \leq   &   \frac{1} { (1-\alpha) \sqrt{1-\delta_{2k}}  } (\| y- A[u^p \otimes (\bigotimes_{j=1}^{\omega} w^{(j)})] \|_2  +  \|\nu'\|_2)     \nonumber   \\
     \leq   &   \frac{  \delta_{2k}    + 2(\omega-1) \delta_{3k}    }  {  1-\alpha  }  \sqrt{\frac{1+ \delta_k}  {1-\delta_{2k}}  }     \|x_S-x^p\|_2 \nonumber \\
   &  +\frac{2\omega-1}{1-\alpha}  \sqrt{\frac{1+ \delta_k}  {1-\delta_{2k}}  }\|A^T \nu'\|_2 + \frac{2}{ (1-\alpha) \sqrt{1-\delta_{2k}} } \|\nu'\|_2.
   \end{align}

\textbf{\emph{Case 2.}} $ \Theta^* > \alpha \|A [  (u^p \otimes  (\bigotimes_{j=1}^{\omega} w^{(j)}) -x_S )_{ X \cup  S} ] \|_2. $
 In this case, by (\ref{EEAA}),   we obtain
 \begin{equation} \label{theta*}  \|[u^p \otimes  (\bigotimes_{j=1}^{\omega} w^{(j)}) -x _S]_{ X \cup  S}\|_2 \leq \frac{\Theta^*} {\alpha \sqrt{1-\delta_{2k}}}.  \end{equation}
So it is sufficient to bound  $\Theta^*.$ The idea is similar to the proof of Lemma  \ref{Lem4455}.
Since  $ (x_S)_{\overline{S\cup X}} =0, $  $\Theta^* $ can be written as
$$ \Theta^*= \| A [( u^p \otimes  (\bigotimes_{j=1}^{\omega} w^{(j)})) _{\overline{X\cup S}}  ]  \|_2   = \| A [((u^p-x_S) \otimes  (\bigotimes_{j=1}^{\omega} w^{(j)}))_{\overline{ X\cup S}} ] \|_2 . $$
Let  $ q $ and $  \kappa $ are integer numbers such that $ \left |\overline{X \cup S} \right |  =  (q-1)  k + \kappa, $ where  $0\leq \kappa < k. $  Let  $(w^{(1)}  )_{\overline{ X\cup S}} $ be decomposed into $k$-sparse vectors as follows:
   $$ (w^{(1)}  )_{\overline{X\cup S}} = (w^{(1)}  )_{S_1} +  \cdots  + (w^{(1)} )_{S_{q-1}}+  (w^{(1)}  )_{S_{q} },  $$
  where  $ {S_1}$ is the index set for  the $k$ largest elements  in  the set $ \left\{(w^{(1)})_i: i \in \overline{X\cup S } \right\},$    and $S_2$ is the index set for the second  $k$ largest   elements in  $ \left\{(w^{(1)})_i: i \in \overline{X\cup S } \right\},$  and so on. $S_q$ is the index set for the remaining $\kappa$ element in this set.   The index sets $S_\ell,   \ell=1, \ldots, q $  are mutually disjoint  and  $  |S_\ell|=k \textrm{ for all } \ell =1, \dots, q-1  $ and $ |S_{q}|= \kappa < k.$  Clearly,  $\overline{X\cup S} = S_1 \cup S_2 \cup  \cdots  \cup S_{q }.  $
  Applying the Lemma \ref{Lem-B2} with $w= w^{(1)}, \tau = k $ and $\Lambda =\overline{X\cup S}$  yields the following inequality:
 \begin{equation}   \label{2222}  \sum_{\ell=1}^{q} \|[ w^{(1)}]_{S_\ell}\|_\infty < 2 . \end{equation}
Define the vector $ z^{(\ell)}: =    [(u^p-x_S) \otimes  (\bigotimes_{j=1}^{\omega} w^{(j)})]_{S_\ell},$  then $$  [(u^p-x_S) \otimes  (\bigotimes_{j=1}^{\omega} w^{(j)})]_{\overline{X\cup S }} = z^{(1)}+ z^{(2)} + \cdots +  z^ {(q)}  . $$
Therefore,
\begin{equation} \label{TTT}  \Theta^*  =    \| A \sum _{ \ell=1}^ {q} z^{(\ell)}  \|_2 \leq \sum _{ \ell =1}^ {q} \|A z^{(\ell )} \|_2 \leq \sqrt{1+\delta_{k}} \sum _{ \ell =1}^{q} \|z^{(\ell)} \|_2,  \end{equation}
where the last inequality follows from the definition of $\delta_k$ and the fact that  every $z^{(\ell)}$ is $k$-sparse. We now estimate the term    $\sum _{ \ell=1}  ^{q} \|v^{(\ell )} \|_2  .$
By (\ref{WWZZ}), we see that
\begin{eqnarray}   \label{6262}   \|z^{(\ell)} \|_2   & =  &  \|[ (u^p-x_S) \otimes (\bigotimes_{j=1}^{\omega} w^{(j)})]_{S_\ell}\|_2  =  \|[((I-A^TA)(x_S-x^p)-A^T\nu') \otimes (\bigotimes_{j=1}^{\omega} w^{(j)})]_{S_\ell} \|_2  \nonumber \\
 &   \leq   & \|[(I-A^TA)(x_S-x^p)]_{S_\ell }  \otimes (\bigotimes_{j=1}^{\omega} w^{(j)})_{S_\ell} \|_2 + \|(A^T\nu')_{S_\ell }  \otimes (\bigotimes_{j=1}^{\omega} w^{(j)})_{S_\ell}  \|_2 \nonumber\\
  & \leq  &    \|(\bigotimes_{j=1}^{\omega} w^{(j)})_{S_\ell}\|_\infty   (\| [(I-A^TA)(x_S-x^p)]_ {S_\ell } \|_2 + \|A^T\nu' \|_2)  \nonumber  \\
  & \leq &  \|(w^{(1)}  )_{S_\ell}\|_\infty   ( \delta_{3k} \|x_S-x^p \|_2 + \|A^T\nu' \|_2)) ,
  \end{eqnarray}  where the last inequality follows from the fact $ 0 \leq w^{(j)}  \leq  \textbf{\textrm{e}}   $ for all $ j=1, \ldots, \omega$ and from Lemma \ref{Lem-Basic} with  $ |S_\ell \cup \textrm{supp} (x_S-  x^p)| \leq 3k. $
 Thus combining (\ref{2222}), (\ref{TTT}) and (\ref{6262}), we obtain
\begin{align*}  \Theta^*     &   \leq     \sqrt{1+\delta_k}  ( \sum _{ \ell =1}  ^{q}  \|[w^{(1)}  ]_{S_\ell}\|_\infty  \left( \delta_{3k} \|x_S-x^p \|_2 + \|A^T\nu' \|_2 \right )  \\
 & \leq 2 \delta_{3k} \sqrt{1+\delta_k}  \|x_S-x^p\|_2 + 2\sqrt{1+\delta_k}\|A^T\nu' \|_2.
\end{align*}
Substituting this into (\ref{theta*}), we get
\begin{equation} \label{Est-08}   \|[u^p \otimes  (\bigotimes_{j=1}^{\omega} w^{(j)}) -x_S ]_{ X \cup  S}\|_2 \leq \frac{2 \delta_{3k} }{\alpha } \sqrt {\frac{ 1+\delta_k  } { 1-\delta_{2k} } }  \|x_S-x^p\|_2 +   \frac{2   }{\alpha } \sqrt {\frac{ 1+\delta_k  } { 1-\delta_{2k} } }  \|A^T \nu'\|_2.  \end{equation}
Thus combining (\ref{Est-02}) for \emph{Case 1} and  (\ref{Est-08}) for \emph{Case 2} yields
\begin{align}  \label{EEDD}  \|[u^p \otimes  (\bigotimes_{j=1}^{\omega} w^{(j)}) -x_S ]_{ X \cup  S}\|_2    & \leq   \max \left\{ \frac{  \delta_{2k}    + 2(\omega-1) \delta_{3k}    }  {  1-\alpha  }     ,  \frac{2 \delta_{3k} }{\alpha }    \right\}  \sqrt{\frac{1+ \delta_k}  {1-\delta_{2k}}  } \|x_S-x^p\|_2  \nonumber  \\
& ~~~ + \max \left\{\frac{2\omega-1 }{1-\alpha}, \frac{2}{ \alpha}\right\} \sqrt{\frac{1+ \delta_k}  {1-\delta_{2k}}  }\|A^T \nu'\|_2 \nonumber \\
& ~~~ + \frac{ 2}{ (1-\alpha) \sqrt{1-\delta_{2k}}} \|\nu'\|_2,
\end{align}
 which holds for any given number $\alpha \in (0,1).$ It is very easy to verify that
 $$  \min_{\alpha \in (0,1)} \max \left\{   \frac{  \delta_{2k}    + 2(\omega-1) \delta_{3k}    }  {  1-\alpha  }     ,  \frac{2 \delta_{3k} }{\alpha }    \right\} =  2 \omega \delta_{3k} +\delta_{2k} . $$  This minimum value attains  at
 \begin{equation} \label{aaaa} \alpha^* =  \frac{ 2 \delta_{3k}  }{  2 \omega \delta_{3k} +\delta_{2k}  }  . \end{equation}
 Combining  (\ref{Bd-01}) and (\ref{EEDD}) produces
 \begin{equation} \label{xxxeee} \|  x _S-x^{\#} \|_2 \leq  \left[ (2 \omega \delta_{3k} +\delta_{2k})  \sqrt {\frac{ 1+\delta_k  } { 1-\delta_{2k} } }  +  \delta_{3k}\right]  \|x^p-x_S\|_2 +c_1 \|A^T\nu'\|_2 + c_2 \|\nu'\|_2,  \end{equation}
where $ c_2  =   \frac{2} { (1-\alpha^*) \sqrt{1-\delta_{2k}}}  $ and  $$ c_1 = \max \left\{\frac{2\omega-1 }{1-\alpha^*}, \frac{2}{ \alpha^*}\right\} \sqrt{\frac{1+ \delta_k}  {1-\delta_{2k}}  } +1  = \frac{2\omega-1 }{1-\alpha^*}  \sqrt{\frac{1+ \delta_k}  {1-\delta_{2k}}  } +1 , $$  where $ \alpha^* $ is given by (\ref{aaaa}) and the second equality above follows from the fact $ \frac{2\omega-1 }{1-\alpha^*} \geq \frac{2}{ \alpha^*}   $ due to the value of $ \alpha^* $ given in (\ref{aaaa}).

 (i) By the structure of the $\textrm{ROT}{\omega},$     $x^{k+1} = x^{\#}. $ Thus the desired result for \textrm{ROT}$\omega$ follows immediately from (\ref{xxxeee}). By noting that $ \delta_k \leq \delta_{2k} \leq \delta_{3k},$  the constant
  $$  \widetilde{\rho}:  =    (2 \omega \delta_{3k} +\delta_{2k})  \sqrt {\frac{ 1+\delta_k  } { 1-\delta_{2k} } }  +  \delta_{3k}       \leq    (2 \omega +1) \delta_{3k} \sqrt {\frac{ 1+\delta_{3k}  } { 1-\delta_{3k} } }  +    \delta_{3k} . $$
  The right-hand side of the above inequality is smaller than 1  provided that $\delta_{3k} < \gamma (\omega), $ where $ \gamma(\omega)     $ is the positive real root in the interval $  (0, 1) $ of the following univariate equation of $\gamma:   $
 $$  g_\omega (\gamma) =   (2 \omega +1)  \gamma \sqrt{ \frac{ 1+\gamma  } { 1-\gamma}  }  +   \gamma  - 1=0.  $$ For a given integer number $ \omega\geq 1, $ the above univariate equation has a unique real root  $\gamma(\omega)$ in the interval $(0,1). $    In fact,   we see that $g_\omega (\gamma)  < 0$ when $\gamma\in (0,1) $ and $ \gamma\to 0,$  and     $ g_\omega (\gamma) >0  $ when $ \gamma\to 1 .$  Also,   the function $ g_\omega (\gamma) $ is strictly increasing over $(0,1). $ This implies that the equation $g_\omega (\gamma) =0 $  has a unique real root in the interval $(0,1).$

 (ii) We now establish the convergence of $\textrm{ROTP}\omega.$ Note that the first step (i.e., the step S1 of the algorithm is the same as that of $\textrm{ROT}\omega.$
  Therefore, the relation (\ref{xxxeee})  remains valid to $\textrm{ROTP}\omega$ which treats $x^\#$ as an intermediate point instead of the next iterate $x^{p+1}. $  Using the point $x^{\#},$ the $\textrm{ROTP}\omega$ algorithm  solve the least-squares problem:
$$ \min_{z} \{ \| y- Az\|_2^2:  ~ \textrm{supp} (z) \subseteq  \textrm{supp} (x^{\#}) \}, $$ to which the solution is set to be  $x^{k+1}.$
By optimality, the vector $x^{p+1}$ must satisfy the relation
$   [A^T (y-  A x^{p+1})]_{\textrm{supp} (x^{\#})} =0    $
which, together with $y=Ax_S+\nu', $  implies that
$$   [(I-A^T A)(x_S-  x^{p+1})]_{\textrm{supp} (x^{\#})} =(x_S-  x^{p+1})_{\textrm{supp} (x^{\#})}+ (\nu')_{\textrm{supp} (x^{\#})},   $$  and hence
\begin{align*} \|(x_S-  x^{p+1})_{\textrm{supp} (x^{\#})}\| &  \leq  \|[(I-A^T A)(x_S-  x^{p+1})]_{\textrm{supp} (x^{\#})}\|_2+ \|(\nu')_{\textrm{supp} (x^{\#})} \|_2 \\
& \leq \delta_{2k}  \| x_S-  x^{p+1}\|_2+ \|\nu'  \|_2.
\end{align*}  The second inequality above  follows from Lemma \ref{Lem-Basic} since $ |\textrm{supp} (x-  x^{p+1})\cup \textrm{supp} (x^{\#}) |\leq 2k. $  Since $\textrm{supp} (x^{p+1}) \subseteq \textrm{supp} (x^{\#})$ which implies $(x^{p+1}-x^{\#})_{\overline{\textrm{supp} (x^{\#})}} =0, $ we then have  that  $$  (x_S-x^{p+1})_{\overline{\textrm{supp} (x^{\#})}} = (x_S-x^{\#}+ x^{\#} - x^{p+1})_{\overline{\textrm{supp} (x^{\#})}}= (x_S-x^{\#})_{ \overline{\textrm{supp} (x^{\#})}}. $$ Therefore,
\begin{eqnarray*}  \|x_S-x^{p+1}\|_2^2 & = &  \|(x_S-x^{p+1})_{\textrm{supp} (x^{\#})} \|_2^2 + \|(x_S-x^{p+1})_{\overline{\textrm{supp} (x^{\#})}}\|_2^2\\
 & \leq  & ( \delta_{2k} \| x_S-  x^{p+1}\|_2+ \|\nu'  \|_2 )^2  +  \|(x_S-x^{\#})_{\overline{ \textrm{supp} (x^{\#})  }}     \|_2^2\\
 & \leq & \delta_{2k}^2 \| x_S-  x^{p+1}\|_2^2 + 2 \delta_{2k} \| x_S-x^{p+1}\|_2 \|\nu'\|_2 + \|\nu'\|_2^2+  \|x_S- x^{\#}  \|_2^2,
\end{eqnarray*}
which can be written as
$$ (1- \delta_{2k}^2) \|x_S-x^{p+1}\|_2^2 - 2 \delta_{2k} \|x_S-x^{p+1}\|_2 \|\nu'  \|_2 - ( \|x_S-x^{\#}   \|_2^2 +\|\nu'  \|_2)  \leq 0. $$
Thus $ \|x_S-x^{p+1}\|_2  $ is smaller than or equal to the largest root of the quadratic equation: $  \phi(t)=:  (1- \delta_{2k}^2) t^2 - 2 \delta_{2k} \|\nu'  \|_2  t -  (\|x_S-x^{\#}   \|_2^2 +\|\nu'  \|_2) = 0. $  This implies  that
\begin{align*} \| x_S-  x^{p+1}\|_2   &  \leq  \frac{1}{\sqrt{1-\delta_{2k}^2}} \|x_S-x^{\#}\|_2 + \frac{1}{1-\delta_{2k}}\|\nu'\|_2 \\
&  \leq  \rho' \|x-x^p \|_2 +\frac{c_1  }{\sqrt{1-\delta_{2k}^2}} \|A^T\nu'\|_2 + \left(\frac{c_2}{\sqrt{1-\delta_{2k}^2 }} +     \frac{1}{1-\delta_{2k}} \right) \|\nu'\|_2.
 \end{align*}
 where the last inequality follows from (\ref{xxxeee}) and the constant $ \rho' $ is given as  $$ \rho'  =\frac{1}{\sqrt{1-\delta_{2k}^2}} \left(   (2 \omega \delta_{3k} +\delta_{2k})  \sqrt {\frac{ 1+\delta_k  } { 1-\delta_{2k} } }  +  \delta_{3k}  \right)  < 1, $$ which is guaranteed  if $\delta_{3k} \leq r^*(\omega),$ where $r^*(\omega)    $ is the real root in $(0,1)$   of the following univariate equation in variable $\gamma :$
 $$   \frac{1}{\sqrt{1-\gamma^2}} \left[   \left( (2 \omega +1)  \gamma \sqrt{ \frac{ 1+ \gamma  } { 1- \gamma}  }  +   \gamma \right)   + \gamma   \right] =1.  $$
 By an analysis similar to (i), it is very easy to verify that the root $\gamma^*(\omega) $ of the above equation in (0,1) is unique. \hfill   $ \Box$

  \vskip 0.08in

 Given a specific integer number   $\omega\geq 1$,  the values of $r(\omega)$ and $r^*(\omega)$ can be immediately obtained.  As a result, the guaranteed performance of ROTP, ROTP2 and ROTP3 (which correspond to the cases $\omega =1,2,3$ respectively) can be immediately obtained from Theorem \ref{Thm-OIHT-02}. For instance, the results for ROTP2 and ROTP3 are   summarized in the corollary below, which is established for the two algorithms for the first time.

\begin{Cor}  \label{finalcor}
  Let $y:  = Ax+\nu $ be the  measurements of $ x$ with measurement errors $\nu. $
   \begin{itemize}

 \item[\emph{(i)}] If   $  \delta_{3k}  \leq  1/7 ,   $ then the sequence $\{x^{p}\} , $ generated by $\textrm{ROTP}2,$  approximates $x_S$ with error
 $$  \| x_S-  x^{p+1}\|_2      \leq  \rho' \|x_S-x^p \|_2 + \gamma_1 \|A^T\nu'\|_2 + \gamma_2 \|\nu'\|_2, $$
  where   \begin{equation} \label{Errf} \rho'  =\frac{1}{\sqrt{1-(\delta_{2k})^2}} \left((\delta_{2k}+ 4\delta_{3k})  \sqrt{ \frac{ 1+ \delta_k } { 1-\delta_{2k}  }}  +  \delta_{3k}\right)  < 1, \end{equation}
 and
 $$ \gamma_1=  \frac{1}{\sqrt{1-\delta_{2k}^2}} \left(1+  \frac{3(\delta_{2k}+ 4 \delta_{3k})} {\delta_{2k} +2\delta_{3k}}  \sqrt{ \frac{ 1+ \delta_k } { 1-\delta_{2k}  }} \right) , $$ $$ \gamma_2 =  \frac{1}{\sqrt{1-\delta_{2k}^2}} \left(\frac{2(\delta_{2k}+ 4 \delta_{3k})} {(\delta_{2k} +2\delta_{3k}) \sqrt{1-\delta_{2k}} } + \frac{1}{1-\delta_{2k}} \right). $$

   \item[\emph{(ii)}]  If   $  \delta_{3k}  \leq  1/9 ,   $ then the sequence $\{x^{p}\} , $ generated by $\textrm{ROTP}3,$  approximates $x_S$   with error
 $$  \| x_S-  x^{p+1}\|_2     \leq   \rho '' \|x_S-x^p \|_2  + \widehat{\gamma}_1  \|A^T\nu'\|_2 +  \widehat{\gamma}_2  \|\nu'\|_2, $$
  where $$  \rho ''   =\frac{1}{\sqrt{1-\delta_{2k}^2}} \left((\delta_{2k} + 6\delta_{3k} )  \sqrt{ \frac{ 1+ \delta_k } { 1-\delta_{2k}  }}  +  \delta_{3k}\right)  < 1, $$ and
   $$ \widehat{\gamma}_1=  \frac{1}{\sqrt{1-\delta_{2k}^2}} \left(1+  \frac{5(\delta_{2k}+ 6 \delta_{3k})} {\delta_{2k} +4\delta_{3k}}  \sqrt{ \frac{ 1+ \delta_k } { 1-\delta_{2k}  }}\right) ,$$ $$ \widehat{\gamma}_2 =  \frac{1}{\sqrt{1-\delta_{2k}^2}} \left(\frac{2(\delta_{2k}+ 6 \delta_{3k})} {(\delta_{2k} +4\delta_{3k}) \sqrt{1-\delta_{2k}} }+ \frac{1}{1-\delta_{2k}}\right). $$
   \end{itemize}
\end{Cor}

The proof of the above corollary is straightforward. In fact,  when $ \omega =2, $ we can verify that  the unique root $r^*(2) $ of the equation (\ref{uni-equ})  in   $(0,1) $ is larger than $1/7.$     For  $ \omega =3 , $ the unique root $r^*(3) $ in  $ (0,1)    $ of (\ref{uni-equ}) is larger than $1/9. $ The corollary    follows from Theorem \ref{Thm-OIHT-02} immediately.
     The  results in \cite{Z19} for  ROT and ROTP (corresponding to $\omega=1$) can be reobtained immediately from  Theorem \ref{Thm-OIHT-02} as well. Similar to Corollary \ref{finalcor}, the first performance results for ROT2 and ROT3 can be obtained immediately from Theorem \ref{Thm-OIHT-02}. Briefly,  the RIP bounds $\delta_{3k} < 1/7$ and  $\delta_{3k} < 1/9$ are the sufficient conditions for the convergence of ROT2 and ROT3, respectively.

The RIP bounds in Theorem \ref{Thm-OIHT-02} for ROT$\omega$ and  ROTP$_\omega$ are the first guaranteed performance criteria developed for the  algorithms based on the concept of optimal $k$-thresholding.     Compared with the sufficient condition $ \delta_{3k} < 1/\sqrt{3}$ for the  guaranteed performance of IHT and  HTP
(\cite{F11}, \cite{FR13}), the sufficient criteria  for  ROT$_\omega$ and  ROTP$_\omega$   are relatively conservative at their current stage. At present, we have only shown that $ \delta_{3k} <1/7 $  is a sufficient condition for the guaranteed performance of ROT2 and ROTP2, and $ \delta_{3k} <1/9 $   for ROT3 and ROTP3.  The appearance of quadratic optimization problem (i.e., compression problem) in ROTP$_\omega$ posts a challenge to the performance analysis of the algorithm,  which is more demanding than that of IHT and HTP.    We believe the  current results  in  Theorem \ref{Thm-OIHT-02} for ROT$\omega$ and ROTP$\omega$  would be improved  if a more convenient and more suitable analysis than the one used in this paper is found.  At the moment, however,  it is not clear whether the RIP bound in Theorem  \ref{Thm-OIHT-02} can be further improved. This is a worthwhile future research topic.

\section{Numerical Performance}
In this section, we discuss the computational complexity of ROTP$_\omega$ and provide the  numerical comparison of this algorithm and  several existing approaches.

\subsection{Computational complexity}  From Theorem \ref{Thm-OIHT-02},  the iterate $ x^p$ generated by   ROTP$_\omega$  approximates the signal with the error
 \begin{equation} \label{EERRFF} \|x^p-x_S\|_2 \leq (\rho')^p \|x^0-x_S\|_2 + (\tau_1 \|A^T \nu'\|_2 + \tau_2 \|\nu'\|_2)/(1-\rho'),  \end{equation}   where $ \rho', \tau_1, \tau_2 $ and $ \nu' $ are  defined in Theorem \ref{Thm-OIHT-02}.   Let $ \varepsilon >0$ be a given tolerance. It follows from (\ref{EERRFF}) that  \begin{equation} \label{rree} \|x^p-x_S\|_2 \leq \varepsilon +  (\tau_1 \|A^T \nu'\|_2 + \tau_2 \|\nu'\|_2)/(1-\rho')   \end{equation}
 provided  $ (\rho')^p \|x^0-x_S\|_2   \leq \varepsilon$ which is guaranteed if
$$ p \geq p^*:= \lceil \log (\varepsilon/\|x^0-x_S\|_2)/\log (\rho') \rceil. $$
The bound (\ref{rree}) indicates that when the signal is  $k$-compressible and its tail $ x_{\overline{S}} $ is small enough, the algorithm   can recover  the significant part of the signal  provided that the measurements are accurate enough, and the algorithm is performed a sufficient number of iterations.
The flops required in one iteration of the algorithm can be estimated as well. To recover a $k$-sparse signal, the governing condition $ \delta_K <1,$ where $k\leq K,$  implies that any $k$ columns of the $m\times n$ matrix are linearly independent, and hence $k \leq m.$  Therefore, we assume $ k \leq m$  in the following complexity analysis of the algorithms, and we only consider the dense measurement matrices for simplicity.

 Obtaining the hard thresholding $ {\cal H}_k$ of a $n$-dimensional vector requires about $ O(n \log k)  $ flops by a sorting approach,  where $ k $ is much smaller than $ n$ in typical compressed sensing scenarios. The projection (least squares) $\min\{\|y-Az\|_2^2: \textrm{supp} (z) \subset \Gamma\} $ with $ |\Gamma | = l \leq m $ is equivalent to solving the normal equation $(A^T_lA_l ) z = A^T_l y , $   which requires about $   m l^2 + l^3/3 \leq \frac{4}{3} m^3 $ flops by using Cholesky decomposition of $ A^T_l A_l. $      The interior-point method in \cite{T88} solving the quadratic  problem (\ref{con-opt})
   requires  $O(n^{3.5} L)$ flops, where $L$ is the size of the problem data encoding in binary.  In each iteration, ROTP$_\omega$ performs one sorting, one projection, and $\omega$ times of the quadratic problem solving. Note that $ \omega$ is a given small integer number independent of $(m,n),$ and  the vector $ u^p= x^p+ A^T(y-Ax^p) $ requires at most $(2m+1)n$ flops.  Thus one  iteration of ROTP$_\omega$  requires about $  O( m^3  + mn + n^{3.5} L)$ flops.

The IHT  only computes the vector $ u^p $  and a hard thresholding of $ u^p$  in each step.  Thus the complexity of one iteration of IHT is about $ O (mn). $
 Since one least-squares problem  is solved in every iteration of HTP, and the HTP requires about  $ O( m^3 + mn)  $ flops in every iteration.  The orthogonal matching pursuit (OMP) needs to perform at least a total of $k$ stages to generate a $k$-sparse vector. It builds up the active set one element at a time, the implementation of OMP can be achieved via updating the Cholesky factorization of the matrix indexed by the active set at leach step. Thus the   $k$ stages of OMP  would take about $  O(m^3+ kmn) $ flops.    The compressive sampling match pursuit (CoSaMP) performs twice of sorting and one projection in every iteration. The total flops needed in one iteration of CoSaMP are about $  O (m^3  + mn).  $ Similarly, the subspace pursuit (SP) needs about $ O(m^3  + mn)  $ flops in every iteration.  The complexity of these algorithms executed  a total of $p$ iterations (except the OMP which performs only $k$ stages) are summarized in the  table below.

\vskip 0.1in

\begin{tabular}{|c||c|}
 \hline
 ~~~~~~Algorithms~~~~~~     & ~~~~~~ Computational Complexity in Dense Matrix Cases ~~~~~~  \\
 \hline
  IHT    &  $   O (pmn)   $    \\
  \hline
 HTP &  $   O(pm^3  + p mn))  $     \\
 \hline
 OMP   &   $   O(m^3+ kmn)       $ \\
 \hline
 CoSaMP   &  $   O( pm^3 +  pmn)     $                         \\
 \hline
 SP     &  $     O(pm^3 +   pmn)     $     \\
 \hline
 ROTP$_\omega$  &  $ O(  pm^3 + pmn  + pn^{3.5} L)    $\\
 \hline
\end{tabular}\\

 The ROTP$_\omega$  algorithm need more computational time  than the others in the above table. However, such extra effort is worthwhile since the emprical results in the next section indicate that  ROTP-type algorithms are usually more reliable  and robust than the above-mentioned existing methods.

\subsection{Comparison with previous methods}

\begin{figure} [htp]
 $ \begin{array}{c}
\includegraphics [width=0.45\textwidth,
totalheight=0.225\textheight] {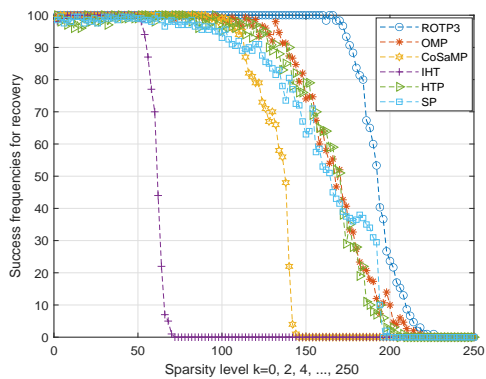} \\
  \textrm{(a) Gaussian measurement matrices}
\end{array} $
 \hfill   $\begin{array}{c}
\includegraphics [width=0.45\textwidth,
totalheight=0.225\textheight] {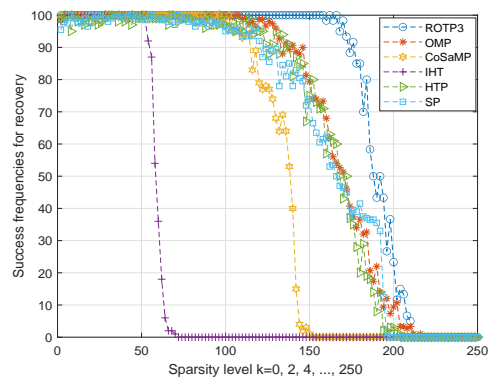} \\
\textrm{(b)  Bernoulli measurement matrices}
\end{array}
$
 \caption{Comparison of the success frequencies of algorithms for signal recovery with inaccurate measurements. For every sparsity level, the success frequency of an algorithm is obtained by  100 random pairs of $ (A, x^*).$  } \label{Fig-A}
\end{figure}

Empirical results show that the traditional IHT and  HTP with stepsize $ \lambda \equiv 1$  is far from efficient for signal reconstruction, and they cannot compete with the relaxed optimal thresholding  methods as shown  in \cite{Z19}.  We now compare the performances of ROTP$_\omega$, OMP, CoSaMP, SP as well as the IHT and HTP with a suitable small stepsize.    In our experiments, all measurement matrices and sparse vectors are randomly generated.
The components of sparse vectors are independent and identically distributed and follow the standard normal distribution, and the positions of nonzero components are uniformly distributed. The size of matrices is set as $400 \times 800,$ and the sparsity level of $k$ of the random vector $ x^* \in \mathbb{R}^{800} $ is ranged from 0 to 220 with stepsize 2, i.e., $ k =0, 2, 4, \dots, 220.$ For every given sparsity level $k$, 100 random pairs $ (A, x^*)$ are realized and used to calculate the success frequency of signal recovery via these algorithms. For every random example $( A, x^*),$  the noisy measurements are given by $ y:= Ax^* + 0.001 h $ where $h$ is a random Gaussian vector. All algorithms in experiment take $x^0=0 $ as the initial point. The recovery criterion is set as $$ \|x^p- x^*\|_2/\|x^*\|_2 \leq 10^{-3}. $$
Two types of measurement matrices were used in simulations:   Gaussian and Bernoulli random matrices.     The results  for  success frequencies of algorithms from Gaussian measurement matrices are summarized in Fig.\ref{Fig-A}(a), and from Bernoulli matrices are given in Fig.\ref{Fig-A}(b).

In our simulations,  the OMP is performed a total of $k$ stages with  $k$ being the sparsity level of the target vector $x^*.$ The parameter $ \omega$ in ROTP$_\omega$ is set to  3, i.e., the specific algorithm ROTP3 is used for comparison. The maximum number of iterations for ROTP3 is set to 40. Due to a low computational cost of  CoSaMP, SP, IHT, and HTP, these algorithms were executed a total of 200 iterations, much more than that of ROTP3 in the experiment.
The efficiency of IHT and HTP relies on the choice of the stepsize $ \lambda,$ and numerical experiments indicate that using  a small stepsize might improve the performance of IHT and HTP in some situations.  However, the performances  of IHT and HTP are quite sensitive to the  choice of stepsize in the sense that a choice of stepsize good in one setting might be inefficient in another environment. Our experiments indicate that the IHT and HTP with  stepsize $\lambda =10^{-3} $ perform  good  for both Gaussian and Bernoulli measurement matrices and for the noisy measurements with errors $\nu= 0.001h,$ where $h$ is a Gaussian noise. So we use this stepsize in IHT and HTP in the performance comparison with other algorithms.

The results in Fig. \ref{Fig-A}   show  that the ROTP3 is an efficient and robust method for signal recovery. It seems more robust than OMP, CoSaMP, SP, IHT and HTP when reconstructing the signal with a relatively high sparsity level $k.$ When $k $ is near $m/2,$ the efficiency of existing algorithms decays very fast, however, the RORP3 can still succeed with a relatively high frequency in signal recovery. It is interesting to observe that   CoSaMP, SP and IHT admits a sharp transition from  high to low success rate at certain level of sparsity in the sense that when the sparsity level $k$ is near some value,  the success rate of recovery for these algorithms might quickly drops to zero.   The decay of the efficiency of OMP and ROTP3  changes  gradually to the variance of sparsity level. This experiment shows that the ROTP3 stands more chance than the  existing methods  to  recover signals with a wider range of sparsity.  Using Bernoulli matrices,  the performance of the algorithm seems slightly different from Gaussian matrices, however, the overall performance in two cases are comparable.

\subsection{Advantage and disadvantage}   The development of ROTP-type algorithms is motivated from the fact that performing the  traditional hard thresholding of an iterate is  independent of the objective function  of  (\ref{L0}).  This weakness might cause the numerical oscillation of the objective function.   The ROTP$\omega$ can be seen as a further development of thresholding technique, whose purpose is to avoid the situation where  a direct use of ${\cal H}_k $ might lead to numerical oscillation. A main feature of the ROTP$\omega$ is to make the thresholding directly  connect to the reduction of objective value, and to ensure ${\cal H}_k$ being applied to a compressed vector, and thus the stability and efficiency of the algorithm are secured. While the data compression  is a quadratic convex optimization which can be efficiently solved by an interior point algorithm,  solving such a data compression problem in each step of ROTP$_\omega$ is clearly time-consuming, compared to the low computational cost of IHT and HTP. This is a disadvantage of ROTP$_\omega$ from a computational point of view. However, the ROTP$_\omega$ would be a good choice in the scenarios when the accuracy and quality of signal recovery are more desired/important than the computational time. As pointed above, the ROTP3 stands more chance than these  existing methods  to  recover signals with a wider range of sparsity.

\section{Conclusions} The newly developed optimal $k$-thresholding algorithms (OT and OTP) can recover $k$-sparse (or $k$-compressible) signals if the restricted isometry constant satisfies $ \delta_k \leq  0.2275 $ when $k$ is even and   $ \delta_{k+1} \leq  0.2275 $ when $k$ is an odd number. Such guaranteed performance conditions governing the sparse  signal recovery are nearly optimal.  Cai and Zhang \cite{CZ13} have proved that $ \delta_k< 1/3 $  is a sufficient condition for the guaranteed recovery of $k$-sparse signals via $\ell_1$-minimization.  A clear   question  is whether the RIP bounds for OT and OTP established in this paper can be improved to $ \delta_k <1/3 $ or $ \delta_{k+1} <1/3 $? Given an integer number $\omega $ (the number of times for data compression in every iteration),   it turns out that the algorithms ROT$\omega $ and ROTP$\omega $ can guarantee to recover the sparse signal if the sensing matrix satisfies the  condition   in Theorem \ref{Thm-OIHT-02}.  As   special cases, the convergence of the  ROTP2 and ROTP3 can be guaranteed under the  bounds $ \delta_{3k} < 1/7  $  and   $ \delta_{3k} < 1/9, $ respectively.  An immediate question is whether these theoretical results can be improved.

\end{document}